\documentclass[11pt]{amsart}
\usepackage{amsmath, amsthm,amssymb}
\usepackage[dvips]{graphicx}
\usepackage{psfrag}

\usepackage[alphabetic]{amsrefs}
\usepackage{times}

\theoremstyle{plain}
\newtheorem{theorem}{Theorem}
\newtheorem{proposition}{Proposition}
\newtheorem{lemma}[theorem]{Lemma}
\newtheorem{corollary}[theorem]{Corollary}

\theoremstyle{definition}

\theoremstyle{remark}
\newtheorem{remark}{Remark}

\long\def\begcom#1\endcom{}

\newcommand{\nc}{\newcommand}

\nc{\sign}{\operatorname{sign}}
\nc{\const}{\operatorname{const\,}}

\nc{\ang}{\sphericalangle}

\nc{\cal}[2]{\ensuremath{\mathcal{#1}^{(#2)}}}
\nc{\tcal}[2]{\ensuremath{\mathcal{\tilde{#1}}^{(#2)}}}
\nc{\hcal}[2]{\ensuremath{\mathcal{\widehat{{#1}}}^{(#2)}}}

\nc{\EG}[1]{\ensuremath{\text{(EG)}_{#1}}}
\nc{\BD}[1]{\ensuremath{\text{(BD)}_{#1}}}
\nc{\HV}[1]{\ensuremath{\text{(HV)}_{#1}}}

\nc{\negfour}{\!\!\!\!} \nc{\ega}[1]{\ensuremath{e^{-\alpha #1}}}
\nc{\egb}[1]{\ensuremath{e^{-2\alpha #1}}}
\nc{\z}[1]{\ensuremath{z_{#1}}}
\nc{\calC}{{\mathcal C}}

\nc{\pfi}[2]{\ensuremath{\partial_{#1}\Phi_{#2}}}
\nc{\pd}[2]{\ensuremath{\partial_{#1}^{#2}}}
\nc{\pdtil}[2]{\ensuremath{\tilde\partial_{#1}^{#2}}}

\begin{document}

\title{Parameter exclusions in H\'enon-like systems}
\author{Stefano Luzzatto}
\address{Mathematics Department,
Imperial College,
London SW7 2AZ, UK}
\email{stefano.luzzatto@imperial.ac.uk}
\urladdr{http://www.ma.ic.ac.uk/\textasciitilde luzzatto}
\author{Marcelo Viana}
\address{Instituto de Matem\'atica Pura e Aplicada, Est. Dona Castorina
110, Rio de Janeiro, Brazil}
\email{viana@impa.br}
\urladdr{http://www.impa.br/\textasciitilde viana}
\thanks{We are most grateful to Sylvain Crovisier and Jean-Christophe
Yoccoz for reading an earlier version and providing very useful
comments. M.V. is partially supported by FAPERJ, Brazil.}
\date{April 12, 2003}

\maketitle

\section{Introduction}

This survey is a presentation of the  arguments
in the proof that H\'enon-like maps
$$
\Phi_a(x,y)=(1-a x^2,0) + R(a,x,y), \quad \|R(a,x,y)\|_{C^3} \le b
$$
have a \emph{strange attractor}, with positive Lebesgue
probability in the parameter $a$, if the perturbation size $b$ is
small enough. We first sketch a \emph{geometric model} of the
strange attractor in this context, emphasising some of its key
geometrical properties, and then focus on the construction and
estimates required to show that this geometric model does indeed
occur for \emph{many parameter values}.

Our ambitious aim is to provide an exposition at one and the same
time intuitive, synthetic, and rigorous. We think of this text as
an introduction and study guide to the original
papers~\cite{BenCar91} and \cite{MorVia93} in which the results
were first proved. We shall concentrate on describing in detail
the overall structure of the argument and the way it breaks down
into its (numerous) constituent sub-arguments, while referring the
reader to the original sources for detailed technical arguments.
Let us begin with some technical and historical remarks aimed at
motivating the problem and placing it in its appropriate
mathematical context.

\subsection{Uniform and non-uniform hyperbolicity}

The arguments which we shall discuss lie at the heart of a certain
branch of dynamics. To formulate its aims and scope we recall
first of all two notions of hyperbolicity: \emph{uniform
hyperbolicity} where hyperbolic estimates are assumed to hold
uniformly at every point of some  set, and \emph{non-uniform
hyperbolicity} which is formulated in terms of asymptotic
hyperbolicity estimates (non-zero Lyapunov exponents) holding only
almost everywhere with respect to some invariant probability
measure. The notion of uniform hyperbolicity was introduced by
Smale (see \cite{Sma67} and references therein) and was central to
a large part of the development the field of Dynamics experienced
through the sixties and the seventies, including the fundamental
work of Anosov~\cite{Ano67} on ergodicity of geodesic flows; the
notion of non-uniform hyperbolicity was formulated by the work of
Pesin \cite{Pes77} and was subsequently much developed by him and
several other mathematicians.

In both cases, one may distinguish two related but distinct
aspects. On the one hand there is the general theory which
\emph{assumes} hyperbolicity and addresses the question of its
geometrical and dynamical \emph{implications} such as the
existence of stable and unstable manifolds, questions of
ergodicity, entropy formulas, statistical properties etc. This
aspect of the theory is well developed in both cases although
results are naturally stronger in the uniformly hyperbolic case.
See the comprehensive texts \cites{Man87,Shu87, Pol93, KatHas94,
AnoSol95, Yoc95a, Via97, Bal00, BarPes01} for details and
extensive bibliographies.

On  the other hand there is the question of constructing and
finding examples and, more generally, of \emph{verifying}
hyperbolicity in specific classes of systems. In this respect, the
difference between uniform and non-uniform hyperbolicity is
striking. Uniformly hyperbolic systems can be  constructed
relatively easily and in principle, and often also in practice, it
is possible to verify the uniform hyperbolicity conditions by
considering only a finite number of iterations of the map. A main
technique for verifying uniform hyperbolicity is the method of
\emph{conefields} which involves checking some open set of
relations on the partial derivatives of the map.

Verifying non-uniform hyperbolicity is generally much more
problematic, partly because this notion is asymptotic in nature,
that is, it depends on the behavior of iterates as time goes to
infinity. Also, non-uniformly hyperbolic systems may contain
tangencies between stable and unstable leaves in which case they
cannot admit complementary stable and unstable continuous
invariant conefields. In fact, invariant objects for this type of
systems tend to live in the measurable category rather than the
topological category. Moreover, there is an a priori
\emph{impasse} related to the fact that the very definition of
non-uniform hyperbolicity requires an invariant measure. Such a
measure is not usually given to begin with and one needs to take
advantage of hyperbolicity features of the system to even prove
that it exists. All in all, we still lack a good understanding of
what makes a dynamical system non-uniformly hyperbolic and it
seems more examples of such systems need yet to be found for such
an understanding to be achieved.

The research which we describe in this paper is at the heart of
ongoing work towards developing a toolbox of concrete conditions
which can play a similar role to that of the conefield conditions
in the uniformly hyperbolic case, implying both the existence of
an invariant measure and the property of non-uniform hyperbolicity
essentially at the same time. The difficulties we mentioned before
are even more significant in the context of H\'enon-like systems
because non-uniformly hyperbolicity cannot be expected to be
persistent in parameter space and thus cannot be checked using an
open set of conditions which only take into account a finite
number of iterations. We shall try to describe here how these
difficulties have been resolved, in a series of spectacular
developments over the last quarter of a century or so.

\subsection{The H\'enon family}
The H\'enon family $H_{a,b}(x,y)=(1-ax^2+by,x)$ was introduced in
the mid-seventies \cite{Hen76} as a simplified model of the
dynamics associated to the Poincar\'e first return map of the
Lorenz system of ordinary differential equations \cite{Lor63} and
as the simplest model of a two-dimensional dynamical system
exhibiting chaotic behavior. The numerical experiments carried out
by H\'enon suggested the presence of a non-periodic attractor for
parameter values $a\approx 1.4$ and $b\approx 0.3$. However,
numerics cannot tell a truly strange attractor from a periodic one
having large period, and rigorous proofs that a strange aperiodic
attractor does exist have proved to be extremely challenging.
H\'enon's original assertion remains unproved to-date for the
parameter range he considered even though remarkable progress has
been made in this direction.

\begin{figure}[htp]
\begin{center}
\includegraphics[height=1.2in]{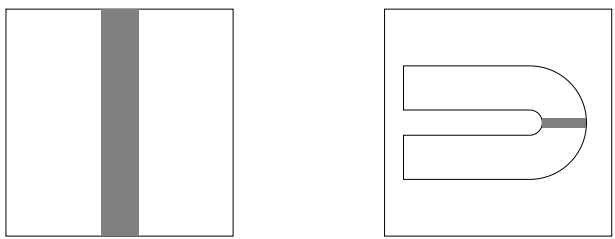}
\caption{\label{f.folding}Folding behavior}
\end{center}
\end{figure}

The distinctive feature of these maps, which makes them a model
for much more general systems, is the occurrence of ``folds"
as described in Figure~\ref{f.folding}: in the
shaded region horizontal/expanding  and vertical/contracting
directions are,
roughly, interchanged. This may give rise to tangencies between
stable and unstable manifolds where
expanding and contracting behavior gets mixed up and
 implies
that if some hyperbolicity is present it will have to be strictly
non-uniform and the dynamics will be structurally unstable.

\subsubsection{The case \( b=0 \)}
In the strongly dissipative limit \( b=0 \), the H\'enon family
reduces to a family of quadratic one-dimensional maps and the fold
reduces to a \emph{critical point}. It is in this context that the
first results appeared. Abundance of aperiodic and  non-uniform
hyperbolic behaviour, was first proved by Jakobson \cite{Jak81}
less than a quarter of a century ago in a paper which pioneered
the parameter exclusion techniques for proving the existence of
dynamical phenomena which occur for nowhere dense positive measure
sets.

The starting point was the formulation of some geometrical
condition which implies non-uniform hyperbolicity; in \cite{Jak81}
this was defined as the existence of an induced map with certain
expansion and distortion properties. A conceptual breakthrough was
the realization that since this condition requires information
about all iterates of the map, it was not reasonable to try to
prove it for a particular given map. Instead one should start with
a family of maps for which some finite number of steps in the
construction of the required induced map can be carried out. One
then tries to take the construction further one step at a time,
and at each step excludes from further consideration those
parameters for which this cannot be done. The problem then reduces
to showing that not all parameters are excluded in the limit and
this is resolved by a probabilistic argument which shows that the
proportion of excluded parameters decreases exponentially fast
with \( n \) implying that the total measure of the exclusions is
relatively small and a positive measure set of parameters remains
for which the constructions can be carried out for all iterations.
All corresponding maps are therefore non-uniformly hyperbolic.

There have been many generalizations of Jakobson's Theorem, using
different geometric conditions to define the notion of a ``good"
parameter, and considering more general families of
one-dimensional smooth maps \cites{BenCar85, Ryc88, MelStr88,
Tsu93, Tsu93a, ThiTreYou94, Cos98, Luz00, HomYou02} as well as
maps with critical points and singularities with unbounded
derivative \cites{PacRovVia98, LuzTuc99, LuzVia00}. Many of these
papers use an intermediate geometric condition formulated in terms
of the properties of the orbits of the critical points which is
sometimes easier to work with than the full induced map. It is
then possible to show by independent arguments that the
appropriate conditions on the orbits of the critical points imply
the existence of an induced map and thus the existence of an
invariant measure and non-uniform hyperbolicity.

\subsubsection{H\'enon-like systems}
Extension of these results to the two-dimensional case requires a
significant amount of new arguments and new ideas. Several issues
will be discussed below when we make a more detailed comparison
with the one-dimensional case. For the moment we just mention the
``conceptual'' problem mentioned above of what a good parameter
looks like. It turns out that it is possible to generalize the
one-dimensional approach mentioned above, formulated in terms of
recurrence properties of the orbits of critical points. However,
even the precise formulation of such a generalization is highly
non-trivial and occupies a central part of the theory. One
outstanding contribution of Benedicks and Carleson~\cite{BenCar91}
was to invent a geometrical structure encompassing tangencies
between stable and unstable leaves, which play the role of
critical points, together with non-uniformly hyperbolic dynamics.
They were then able to generalize the parameter exclusion argument
to conclude that this structure does occur in the H\'enon family
$H_{a,b}$\,, for a positive Lebesgue measure set of parameters
$(a,b)$ with $b \approx 0$.

Shortly afterwards, \cite{MorVia93} extended Benedicks and
Carleson's approach to general H\'enon-like families, thus freeing
the arguments from any dependence on the explicit expression of
the H\'enon maps, and also established the connection between
these systems and general bifurcation mechanisms such as
homoclinic tangencies. Moreover, \cite{Via93} extended the
conclusions of \cite{MorVia93} to arbitrary dimension. The ergodic
theory of H\'enon-like systems was then developed, including the
existence of a Sinai-Ruelle-Bowen measure~\cite {BenYou93} (in
particular proving that the attractors of
\cites{BenCar91,MorVia93} are indeed non-uniformly hyperbolic in
the standard sense), exponential mixing \cites{BenYou00,You98},
and the basin property \cite{BenVia01}. Moreover,
\cite{DiaRocVia96} extended \cite{MorVia93} to the, more global,
strange attractors arising from saddle-node cycles and, in doing
so, observed that the original approach applies to perturbations
of very general families of uni- or multimodal maps in one
dimension, besides the quadratic family. Recently, \cite{WanYou01}
showed that the whole theory extends to such a generality, and
also isolated a small set of conditions under which it works.
Similar ideas have also been applied in related contexts in
\cites{Cos98, PumRod01, PalYoc01, WanYou02}.

\subsection{General remarks and overview of the paper}

One key point in the construction in \cite{BenCar91} is the notion
of \emph{dynamically defined critical point}, a highly non-trivial
generalization of the notion of critical point in the
one-dimensional context, and the associated notion of dynamically
defined \emph{finite time approximation} to a critical point. The
definition of a good parameter is formulated in terms of the
existence of a suitable set of such critical points satisfying
certain hyperbolicity conditions along their forward orbits.
However the very existence of the critical points is tied to their
satisfying such hyperbolicity properties and thus to the parameter
being a good parameter, and we are faced with another
\emph{impasse} analogous to the one discussed above.

The solution lies in the observation that a set of finite time
\emph{approximations} to these critical points can be defined for
all parameter values in some sufficiently small parameter
interval. One can then set up an inductive argument where a
certain condition satisfied by the critical approximations implies
that the approximations can be refined to a better approximation.
Parameters for which the condition is not satisfied are excluded
from further consideration. Then, as in the one dimensional case,
one has to estimate the size of the exclusions at each step to
conclude that there is a substantial set of parameters for which
all critical approximations always satisfy the required condition
and in particular converge to a ``true'' critical set which also
satisfies these conditions.

The overall argument is set up as an induction which is quite
involved, and the exposition in the original papers is
occasionally terse, especially when describing the parameter
exclusions procedure. More explanations on some important points
have been provided subsequently, for instance in
\cite{PacRovVia98}, where the handling of infinitely many critical
points was formalized in detail, in a one-dimensional set-up.
However, it has been suggested that it would be useful to have in
a single text a conceptual reader-friendly survey of the whole
procedure with particular emphasis on parameter exclusions. The
present text is a response to that suggestion.

In Section \ref{s.induction} we briefly outline the main
geometrical properties of the attractor for the Benedicks-Carleson
``good'' parameter values, including the definition of dynamically
defined critical points. This corresponds to formulating precisely
the conditions which determine the parameters which will be
excluded at each step of the parameter exclusion argument. In the
remainder of the paper we discuss the second stage: proving that
the set of good parameters has positive Lebesgue measure.

Another text, with a similar goal, has been written at about the
same time by Benedicks, Carleson \cite{BenCar02}, and another
presentation of parameter exclusions is contained in
\cite{WanYou01} in a related more general setting. Our
presentation is based on the original papers \cite{BenCar91} and
\cite{MorVia93}, although we present here a new (previously
unpublished) formalization of the arguments by introducing the
notion of an \emph{extended parameter space} to keep track of the
combinatorics of each individual critical point approximation at
each stage and to make more explicit the effect of exclusions due
to multiple critical points. This formalism was developed as part
of ongoing joint work on Lorenz-like attractors \cite{LuzVia} and
was first announced in \cite{Luz98}.

\section{Geometrical structure in dynamical space}\label{s.induction}

In this section we review the basic geometric properties of a
``good'' parameter value, and introduce the notation and definitions
required to set up the parameter exclusion argument.
When this is not a source of confusion we
will consider the parameter $a$ to be  fixed and will
not mention it  explicitly.

\subsection{The one-dimensional case}\label{ss.hyperbolicity}
The argument in the one-dimensional case $\phi_a:x\mapsto 1-ax^2$
breaks down into three basic steps.

\subsubsection{Uniform expansion outside a critical neighbourhood}

The first step is a manifestation of the general principle in one
dimensional dynamics, proved by Ma\~n\'e~\cite{Man85}, according
to which orbits behave in a uniformly hyperbolic fashion as long
as they remain outside a neighborhood of the critical points and
the periodic attractors. More specifically, we use

\begin{proposition}\label{1dimunifexp}
    There exists a constant \( \lambda>1 \) such that
    for every \( \delta>0 \) there exists \( a_{0}(\delta) <2 \)
    such that for every \( a_{0}\leq a \leq 2 \),
    the dynamics of \( \phi_{a} \) outside a \( \delta \)-neighbourhood of the
    critical point $c=0$ is uniformly expanding with expansion rate \(
    \lambda \).
\end{proposition}
Thus by choosing a parameter interval \( \Omega \) sufficiently
close to \( a=2 \) we can work with maps which satisfy uniform
expansion estimates, uniformly also in the parameter, outside any
\emph{arbitrarily small} neighbourhood of the critical point with
an expansion coefficient \( \lambda \) which does not depend on
the size of the neighbourhood. This fact is crucial to the whole
sequel of the arguments.

\subsubsection{Bounded recurrence and non-uniform expansivity}
Once a constant \( \delta>0 \), the corresponding critical
neighbourhood \( \Delta \), and
a suitable parameter interval \( \Omega \) have been fixed, we
define a good parameter \( a \) by the recurrence condition
\begin{equation}\tag{\( * \)}
\sum_{\substack{
1\leq j\leq n \\
c_{j}\in\Delta
}}
\log |c_j-c|^{-1} \le\alpha n \quad\text{ for all } n\ge 1,
\end{equation}
where \( c_{j}= \phi_a(c)\) are the iterates of the critical
point, and $\alpha>0$ is some small constant. We remark that this
is different, but essentially equivalent to, the \emph{basic
assumption} and the \emph{free period assumption} taken together,
which are the conditions originally formulated in \cites{BenCar91,
MorVia93}. It is similar to the condition of \cite{Tsu93b} and has
proved particularly useful in \cites{LuzTuc99, LuzVia00, LuzVia}
where the presence of a discontinuity set requires an additional
bounded recurrence condition which remarkably takes exactly the
same form. Moreover, a straightforward calculation using the
expansion estimates of Proposition \ref{1dimunifexp}, see
\cite{Luz00}, shows that under this condition the critical orbit
exhibits exponential growth of the derivative:
\begin{equation}\tag{{\rm EG}}
|D\phi^n(\phi(c))| \ge e^{\kappa n} \quad\text{ for all } n\ge 1,
\end{equation}
for some constant $\kappa>0$. By \cites{ColEck83,NowStr88}
condition \( (EG) \) implies that the corresponding map is
non-uniformly hyperbolic.

\subsubsection{Parameter exclusions}
Thus the problem has been reduced to showing that
many parameters in \( \Omega \) satisfy the bounded recurrence
condition \( (*) \). This is essentially a consequence of Proposition
\ref{1dimunifexp} and the observation that the uniform expansion
estimates given there transfer to expansion estimates for the
derivatives with respect to the parameter. Thus the images of the
critical orbit for different parameter values tend to be more and
more ``randomly'' distributed and the probability of them falling very
close to the critical point gets smaller and smaller. Thus the
probability of satisfying the bounded recurrence conditions is
positive even over all iterates.  In section \ref{ss.one-dimensional}
we sketch the combinatorial construction and the estimates required
to formalize this strategy.  This is also a special case of the
strategy applied to the two-dimensional case which will be discussed
in some detail.

\subsection{The two-dimensional case}

In the two dimensional situation we can also break down
the overall argument into three steps as above, although each one is
significantly more involved. In particular, the very formulation  of
the recurrence condition \( (*) \) requires substantial work and we
concentrate on this issue here, leaving the issues related to the
exclusion of parameters to the later sections.

\subsubsection{Uniform hyperbolicity outside a critical neighbourhood}
In two dimensions we define the critical neighbourhood  \( \Delta \)
as a small   vertical strip of width \( 2\delta \).

\begin{proposition}\label{2dimunifexp}
There exists a constant \( \lambda>1 \) such that for every \(
\delta>0 \) there exists \( b_{0}(\delta)>0 \) and \(
a_{0}(\delta) <2 \) such that for every \( 0\leq b \leq b_{0} \)
and \( a_{0}\leq a \leq 2 \), the dynamics of \( \Phi_{a, b} \)
outside a vertical strip $ \Delta $ around $x=0$ is uniformly
hyperbolic with expansion rate \( \lambda \) and contraction rate
\( b \).
\end{proposition}

The proof of this proposition relies on the fact that \( b \) is
small and thus \( \Phi_{a,b} \) is close to the one dimensional
family of maps for which the estimates of Proposition
\ref{1dimunifexp} hold. One other place in which the strong
dissipativeness assumption \( b\approx 0 \) is used is for
estimating the cardinality of the critical set at each step $n$,
as we shall see.

\subsubsection{Bounded recurrence and non-uniform hyperbolicity}

We now suppose that the constants
\( \delta, a_{0}, b_{0} \) are fixed and that for some \( 0< b \leq
b_{0} \) we have chosen an  interval \( \Omega\subset (a_{0}, 2) \) of
\( a \)-parameters. We want to formulate some
condition with which to characterize the good parameters in \( \Omega
\). Our aim is to remain as close as possible to the
one-dimensional formulation, and to
identify a \emph{critical set} \( \mathcal
C\subset \Delta \) containing an infinite number of \emph{critical
points} such that each point satisfies the bounded recurrence
condition
\begin{equation}\tag{\( * \)}
\sum_{\substack{
1\leq j\leq n \\
c_{j}\in\Delta
}}
\log |c_j-\mathcal C|^{-1} \le\alpha n \quad\text{ for all } n\ge 1,
\end{equation}
for some \( \alpha>0 \)  sufficiently small. Here the distance
\( |c_j-\mathcal C| \) does not refer exactly to the standard
Hausdorff distance between the point \( c_{j} \) and the set \(
\mathcal C \) but to the distance between \( c_{j} \) and some
particular point of \( \mathcal C \) which is chosen by a procedure
to be discussed below.

A not-so-straightforward calculation (which is in fact a large
part of the proof of Theorem \ref{positive exponents} below) shows
that this bounded recurrence condition implies the two dimensional
analogue of the exponential growth condition:
\begin{equation}\tag{{\rm EG}}
|D\phi^n_{\phi(c)}w| \ge e^{\kappa n} \quad\text{ for all } n\ge 1,
\end{equation}
for some constant $\kappa>0$ and for a horizontal or ``almost
horizontal'' vector \( w \) and for any critical point \(
c\in\mathcal C \). We shall not discuss here why \( (*) \)
(together with the uniform hyperbolicity conditions outside \(
\Delta \)) is also  sufficient to guarantee the global non-uniform
hyperbolicity of the corresponding map, and refer to the papers
\cites{BenYou93, BenYou00, BenVia01, HolLuz} for the construction
of the Sinai-Ruelle-Bowen measure under these, or other
essentially equivalent, conditions. We also postpone the
discussion of the verification that this condition is satisfied by
many parameters in \( \Omega \) to the following sections.
Instead, in the remaining parts of this section we focus on the
problem of the definition of the critical set \( \mathcal C \).

\subsubsection{Dynamically defined critical points}
Since \( \Phi \) is a diffeomorphism its Jacobian never vanishes
and thus there are no a priori given critical points as in the
one-dimensional case. However something ``bad'' does happen in the
critical region because the uniformly hyperbolic estimates outside
\( \Delta \) cannot be extended to \( \Delta \). Geometrically
this is due to the folds described above, which are reflected at
the level of the differential by the fact that (roughly)
horizontal vectors get mapped to (roughly) vertical vectors.
Dynamically this is problematic because the (roughly) horizontal
direction is expanding while the (roughly) vertical direction is
strongly contracting. Thus any expansion gained over several
iterates may be lost during the iterates following a return to \(
\Delta \). It is necessary to have a finer control over the way in
which vectors rotate in order to show that after some bounded time
and some bounded contraction, they return to a (roughly)
horizontal direction and start expanding again.

One can hope to characterize geometrically as \emph{critical
points} those points on which the fold has the most dramatic
effect, i.e. those for which the almost horizontal vector which is
\emph{precisely} in the expanding direction, i.e. that vector
which is tangent to an unstable manifold, gets mapped to the
almost vertical vector which is precisely in a contracting
direction, i.e. tangent to a stable manifold. This turns out
indeed to be the case and the set of critical points \( \mathcal C
\) is formed by a set of points of tangential intersection between
some stable and some unstable manifolds. However, as mentioned
above, such manifolds cannot be assumed to exist for all parameter
values, and thus  the construction
 requires an inductive approximation argument by
\emph{finite time critical points} which are also tangencies between
pieces of unstable manifold and some \emph{finite time stable leaves}
to be described below.
In the following sections we shall explain in more detail the
 local geometry associated to critical points and their approximations.
 For the moment we clarify the formal structure of the induction.

 \subsubsection{The induction}
We start by defining a critical set \( \mathcal C^{(0)} \) and
then suppose inductively that a set $\mathcal C^{(k)}$ of critical
points of order $k$ is defined for $1\leq k\leq n-1$, such that
each critical point satisfies certain hyperbolicity conditions
\EG{n-1} which are finite time versions of condition \( (EG) \)
given above, together with a condition \BD{n-1} of bounded
distortion in a neighbourhood to be stated below. The existence of
the set \( \cal C{n-1} \) allows us to \emph{state} a condition
$(*)_{n-1}$ on the recurrence of points of \( \cal C{n-1} \) to
the set itself. This is a finite time version of condition \( (*)
\) given above, with \( \cal C{n-1} \) replacing \( \mathcal C \).

The main inductive step then consists of showing that if all
points of \( \cal C{n-1} \) satisfy this recurrence condition,
then conditions \EG n and \BD n hold in a neighborhood of
$\mathcal C^{(n-1)}$. Now the fact that these conditions hold is
enough to allow us to \emph{define} a new critical set $\mathcal
C^{(n)}$ close enough in the Hausdorff metric to $\mathcal
C^{(n-1)}$ so that its points also automatically satisfy
$(*)_{n-1} $ and \EG n and \BD n. This completes the inductive
step. If all points of $\mathcal C^{(n)}$ satisfy $(*)_{n}$ the
argument can be repeated to obtain a critical set $\mathcal
C^{(n+1)}$ and so on. The sets $\mathcal C^{(n)}$ eventually
converge to a critical set $\mathcal C$ which consists of
tangencies between stable and unstable leaves. We summarize this
reasoning in the following

\begin{theorem}
\label{positive exponents}
Suppose that for some \( a\in\Omega\)
a finite critical set \( \cal C{n-1} \subset \Delta \) has been defined.
\begin{enumerate}
\item
If \( \cal C{n-1} \) satisfies \( (*)_{n-1} \) then it satisfies \EG n
and \BD n;
\item
If \( \cal C{n-1} \) satisfies \( (*)_{n-1} \), \EG n and \BD n, then  a
finite set
\( \cal Cn \) can be defined  whose elements
are critical points of order \( n \)
and  satisfy \( (*)_{n-1} \), \EG n and \BD n.
Moreover \( \cal Cn \) and \( \cal C{n-1} \) are exponentially close in \( n \) in the Hausdorff
sense.
\end{enumerate}
In particular, if \( (*)_{n} \) continues to hold
 for increasing values of \( n \) the set of critical approximations
 converges to a set \( \mathcal C \) of \emph{true} critical points
  satisfying \( (*)_{n} \), \EG n and \BD n for all
 \( n \).
\end{theorem}

This result tells us that the bounded recurrence conditions \(
(*)_{n} \) are exactly the conditions we need to define a good
parameter, and
 allows us to focus the parameter exclusion argument on the
recurrence of the critical approximations at each stage \( n \).
The framework is henceforth similar to the one-dimensional case
apart from the additional complications coming from the
requirement to prove the inductive step and the fact that the
exclusions need to be made with respect to each critical point.

\subsubsection{Why do we need a critical \emph{set} ?}
The reason one needs a whole critical set, and not just a single
critical point, is the way iterates hitting the critical region
$\Delta$ are compensated for in order to recover exponential
growth. Whenever a point $z\in\cal C{n}$ returns to $\Delta$ at
time $n$ one looks for some point $\zeta\in\cal C{n}$ close to
$\Phi^n(z)$, and transmits information about hyperbolicity on the
first iterates of $\zeta$, inductively, to the stretch of orbit
$z$ that follows the return. This works out well if the two points
$\zeta$ and $\Phi^n(z)$ are {\em in tangential position}, that is,
contained in the same almost horizontal curve. For this, in
general, $\zeta$ must be different from $z$. This step forces the
critical sets $\cal Cn$ to be fairly large, indeed, their
cardinality has to go to infinity as $n\to\infty$. Fortunately, as
we are going to see, one can do with a sequence of critical sets
whose cardinality grows slowly enough, as long as one supposes
that $b$ is small. See Section~\ref{cardinality}. In the sequel we
define more formally the notion of critical point and sketch the
argument in the proof of Theorem \ref{positive exponents}.

\subsubsection{Constants and notation}\label{sss.constants}

First we introduce some notation which will be used extensively
below. Given a point $\xi_{0}$ and a  vector $w_0$ we denote
$\xi_j=\Phi^{j}(\xi_{0})$ and $w_j=D\Phi^j(\xi_0)w_0$ for all
$j\ge0$. The vector \( w_{0} \) will be assumed to have slope \(
\leq 1/10 \) unless we explicitly mention otherwise. We fix $\log
2 > \kappa \gg \alpha \gg \delta>0$. These constants have the
following meaning:
\begin{itemize}
\item $\kappa$ is a lower bound for the hyperbolicity of the
two-dimensional map in condition \EG n;
\item $\alpha$ is used in formulating the recurrence condition $ (*)_{n}
$, as well as in defining the notion of binding;
\item $\delta$ defines the width of the critical neighborhood $\Delta$.
\end{itemize}
The parameter interval $\Omega$ is chosen close enough to $a=2$
depending on $\kappa$, $\alpha$, $\delta$. The perturbation size
$b$ is taken to be small, depending on all the previous choices.
A few ancillary constants appear in the course of the arguments,
related to the previous ones. $K=5$ is an upper bound for the
$C^3$ norm of our maps. A small $\rho>0$ e.g. $\rho=(10K^2)^{-2}$
is used to describe the radius of an admissible segment of
unstable manifold around every critical point. We use $\tau>0$ in
the treatment of the recurrence condition. It is chosen in
Section~\ref{sss.verifying}, much smaller than $\alpha$ and
independent of $\delta$. Constants $\kappa_1, \ldots, \kappa_4>0$
depending only on $\kappa$ describe expansion during binding
periods. And $\theta=C/|\log b|$ is used when bounding the number
of critical points, where $C>0$ is some large constant e.g.
$C=10\log|\rho|$. Notice that $\theta\to 0$ as $b\to 0$.

\subsection{Hyperbolic coordinates}
The definition of finite time critical point is based on the
notion of hyperbolic coordinates which we discuss in this section.

\subsubsection{Non-conformal linear maps}
Suppose that the derivative map \( D\Phi^{k}_{\xi_{0}} \) at some
point \( \xi_{0} \) is non-conformal (a very mild kind of
hyperbolicity). Then there are well defined \emph{orthogonal}
subspaces \( E^{(k)}(\xi_{0}) \) and \( F^{(k)}(\xi_{0}) \) of the
tangent space, for which vectors are \emph{most contracted} and
\emph{most expanded} respectively by \( D\Phi^{k}_{\xi_{0}} \).
This follows by the elementary observation from linear algebra
that a linear map \( L \) which sends the unit circle \( \mathcal
S^{1} \) to an ellipse \( L (\mathcal S^{1}) \neq \mathcal S^{1}
\) defines two \emph{orthogonal} vectors \( e \) and \( f \) whose
images map to the minor and major axis of the ellipse
respectively.  The directions \( E^{(k)} \) and \( F^{(k)} \) can
in principle be obtained explicitly as solutions to the
differential equation \({d}|D\Phi^{k}_{\xi_{0}}\cdot (\sin\theta,
\cos\theta)|/{d\theta}=0\) which gives
\begin{equation}\label{c4}
\tan 2\theta  =
\frac{2 (\pfi x1^{k}\pfi y1^{k} +\pfi x2^{k}\pfi y2^{k})}
{(\pfi x1^{k})^2+(\pfi x2^{k})^2 - (\pfi y1^{k})^2 -(\pfi y2^{k})^2}.
\end{equation}
This shows  that the direction fields given by \( E^{(k)} \)
and \( F^{(k)} \) depend smoothly on the base
point and extend to some neighbourhood
of \( \xi_{0} \) on which the derivative continues to satisfy the
required non-conformality.
Therefore  they can be integrated to
give two smooth \emph{orthogonal} foliations
\( \mathcal E^{(k)} \) and \( \mathcal F^{(k)} \).
The individual leaves of these foliations are the
natural \emph{finite time version} of classical local stable and
unstable manifolds. Indeed they are canonically defined precisely
by the property that they are the most contracted and most expanded
respectively for a certain finite number of iterations.
The estimates to be developed below will show that in certain
situations the stable leaves \( \mathcal E^{(k)}(\xi_{0}) \)
converge as \( k\to\infty \) to the classical local stable manifold \(
W_{\varepsilon}^{s}(\xi_{0}) \).

The notions of most contracted directions and most contracted
integral curves play a central role in the original papers
\cites{BenCar91, MorVia93}, although they are not exploited as
systematically as in here. Our formalism was developed in the
context of Lorenz-like systems \cites{LuzVia, HolLuz} and leads to
a significant simplification of several steps of the argument, in
particular it plays an important role in allowing us to formulate
the induction of Theorem \ref{positive exponents} in such a
straightforward way. The idea of approximating the classical local
stable manifold by finite time local stable manifolds has been
further refined in \cite{HolLuz03} where it forms the basis of a
new approach to the local stable manifold theorem in more
classical contexts.

\subsubsection{Notation}
Before explaining how
these foliations are used to define the notion of critical point, we
introduce some more notation.
We let
\(
H^{(k)} (\xi_{0}) = \{ F^{(k)}(\xi_{0}), E^{(k)}(\xi_{0}) \}
\)
denote the coordinate system in the tangent space at \( \xi_{0} \)
determined by the directions \( F^{(k)} \) and \( E^{(k)} \) and by
\(
\mathcal H^{(k)} = \{\mathcal F^{(k)}, \mathcal E^{(k)}\}
\)
the family of such coordinate systems in the neighbourhood in which
they are defined.
We also let
\(f^{(k)} (\xi_{0}) \) and \(e^{(k)}(\xi_{0}) \)
denote unit vectors in the directions \( F^{(k)}(\xi_{0}) \) and
\( E^{(k)}(\xi_{0}) \) respectively.  For \( k=0 \) we let \( F^{(0)} \)
and \( E^{(0)} \) denote the horizontal and vertical direction
respectively and thus \( \mathcal F^{(0)} \) and \( \mathcal E^{(0)}
\) are horizontal and vertical foliations respectively.
Notice that \( \mathcal H^{(k)} \) can be thought of  as living in the
tangent bundle as a family of coordinate systems, or in the phase
space as a foliation; we will not distinguish formally between these two
interpretations.
For all the objects defined
above we use a subscript \( j \)
to denote their images under the map \( \Phi^{j} \), or the
differential map \( D\Phi^{j} \) as appropriate. In particular we let
\( e^{(k)}_{j} = D\Phi^{j} (e^{(k)}),  f^{(k)}_{j} = D\Phi^{j}
(f^{(k)})\)
and
\( \mathcal H^{(k)}_{j} =  \Phi^{j} \mathcal H^{(k)}\).
 Notice moreover, that
\( \mathcal H^{(k)}_{k} \) is also an orthogonal system of
coordinates, whereas
\( \mathcal H^{(k)}_{j} \) is \emph{not} orthogonal in general for \( j\neq k
\).  Notice that the differential map
 \( D\Phi^{k}_{\xi_{0}} \), expressed as a
matrix with respect to the
hyperbolic coordinates \( \cal Hk \) and \( \cal Hk_{k} \), has the
diagonal form
 \[
D\Phi^{k}_{\xi_{0}}= \begin{pmatrix}
\|f^{(k)}_{k}(\xi_{0})\| & 0\\ 0& \|e^{(k)}_{k}(\xi_{0}) \|
\end{pmatrix}.
 \]
Finally, for \( j, k\geq 0 \), we consider  the angle between the
leaves of \( \mathcal H^{(k)} \) and \( \mathcal H^{(k+1)} \) at some
point \( \xi_{0} \) at which both foliations are defined, and the
corresponding angle between the images: \[
\theta^{k}=\theta^{(k)}(\xi_{0}) = \ang (e^{(k)}(\xi_{0}),
e^{(k+1)}(\xi_{0}) ) \ \ \text{ and } \ \
\theta^{(k)}_{j}=D\Phi^{j}_{\xi_{0}} \theta^{(k)}(\xi_{0}) ;
 \]
as well as the derivatives of these angles with respect to the base
point \( \xi_{0} \):
\[
D\theta^{(k)}= D_{\xi_{0}}\theta^{(k)}(\xi_{0})
\ \ \text{ and } \ \
D\theta^{(k)}_{j} = D_{\xi_{0}}\theta^{(k)}_{j}(\xi_{0})
\]
\subsubsection{Convergence of hyperbolic coordinates}
For \( k=1 \) and \( \xi_{0}\notin\Delta \), relation \eqref{c4}
implies
\begin{equation}\label{anglezero}
    |\theta^{(0)}|= \mathcal O(b)
    \end{equation}
    Thus by taking \( b \) small we can guarantee that the stable
    and unstable foliations \( \mathcal E^{(1)}, \mathcal F^{(1)} \)
    are arbitrarily close to the vertical and horizontal foliations
    respectively.
It turns out that the angle between successive contractive directions
is related to the hyperbolicity along the orbit in question,
and we get a quite general
estimate which says that as long as the inductive assumption
\EG k
continues to be satisfied,we have
\begin{equation}\label{anglek}
 |\theta^{(k)}|= \mathcal O(b^{k}).
\end{equation}
In particular, in the limit they converge to a well defined
direction which is contracted by all forward iterates. This
convergence is in the $C^1$ norm (even $C^r$ for any fixed $r$) if
$b$ is sufficiently small.

 \subsection{Critical points}\label{ss.properties}
We are now ready to define
the notion of \emph{critical point of order}
$k$, generally denoted by $ z^{(k)} $. The definition will be given
inductively. The set of such points will
be denoted $\cal Ck$ and is contained in
the critical neighborhood $\Delta$ defined by
\[
\Delta=\{(x,y): |x| \leq \delta \}.
\]
All critical points are on the global unstable manifold $W^u(P)$
of the hyperbolic fixed point $P\approx (1/2,0)$; notice that \(
W^u(P) \) has many folds and \( W^u(P)\cap\Delta \) has infinitely
many connected components. Fix a compact admissible neighborhood
$W$ of $P$ inside $W^u(P)$ with length $\lesssim 2$ and extending
to the left of $P$ across the critical region $\Delta$. By
definition, the intersection of $W$ with the vertical line
$\{x=0\}$ is the unique critical point $z^{(0)}$ of order zero.

\subsubsection{First step of the induction}
Now consider the curve \( W_{0}=\Phi (W\cap \Delta). \) As
described above, the map \( \Phi \) gives rise to a \emph{fold}
precisely in \( \Delta \) and therefore \( W_{0} \) is
\emph{folded} horizontally. The quadratic nature of \( \Phi \)
guarantees that it is in fact a quadratic parabola (positive
curvature) laying on its side. Notice moreover that \(
W_{0}\cap\Delta = \emptyset \) and therefore, contractive
directions \( e^{(1)}_{\xi_{0}} \) of order 1 are defined at each
point \( \xi_{0}\in W_{0}\). The smoothness of these directions,
the fact that they are essentially vertical, see
\eqref{anglezero}, and the fact that \( W_{0}\) is quadratic,
guarantee that there must be a point \( z^{(1)}_{0} \in W_{0}\)
which is tangent to a contracting leaf of the foliation \(
\mathcal E^{(1)} \). A bit more work shows that the leaves of the
contracting foliation \( \mathcal E^{(1)} \) have  small curvature
which further implies that there can be at most one point of
tangency and that this tangency is quadratic.

We define \( z^{(1)}=\Phi^{-1}(z^{(1)}_{0}) \) as the unique
critical point of order 1, i.e. the unique element of the set \(
\mathcal C^{(1)} \), and \( z_{0}^{(1)} \) as the corresponding
critical value. Notice that taking \( b \) small implies that the
``tip'' of the parabola is close to \( z^{(0)}_{0} \) and that the
stable leaves are almost vertical. Therefore the point of tangency
\( z_{0}^{(1)} \) is close to \( z_{0}^{(0)} \) and
 the distance between the
critical points \( z^{(0)} \) and \( z^{(1)} \) is \( \mathcal O(b)
\).
This constitutes the first step in the inductive definition of the
critical set.

The characteristic feature of a critical point $z^{(k)}$ of order
$k$ will be that the unstable manifold is tangent to the stable
foliation of order $k$ at the critical value
$z_0^{(k)}=\Phi(z^{(k)})$. More formally, we assume that for each
$ k= 1, \ldots, n-1 $, the critical set $\cal Ck$ contains points
$z^{(k)}$ with the following properties.

\subsubsection{Generation of critical points}
We introduce the notion of the \emph{generation} of a critical
point which is quite different from the notion of the \emph{order}
of the critical point. We say that the critical point \( k \) is
of generation \( g \geq 1 \) if it belongs to \(
\Phi^{g}(W)\setminus \Phi^{g-1}(W) \), where \( W \) is the
component of \( W^u(P) \) defined above, of length \( \lesssim 2
\) containing the fixed point \( P \) and crossing \( \Delta \)
completely. By convention we say that a critical point is of
generation \( 0 \) if it belongs to \( W \). The critical points
\( z^{(0)} \) and \( z^{(1)} \) defined above, are critical points
of generation \( 0 \) and so are all their refinements \( z^{(k)}
\), \( k\geq 1 \) to be defined below. As part of the construction
we impose the condition that critical points of order \( k \) must
be of generation \( \leq \theta k \), where \( \theta \approx
1/\log b^{-1} \). In particular the only admissible critical
points of order \( \lesssim \log b^{-1} \) are those of generation
0.

\subsubsection{Admissible segments}\label{sss.admissible}

The neighborhood $\omega$ of radius $\rho^{\theta k}$ around
$z^{(k)}$ inside the unstable manifold $W^u(P)$ is an admissible
curve contained in $\Phi^{\theta k}(W)\cap\Delta$; we say that a
curve is almost horizontal, or {\em admissible\/}, if it is a
graph $\{x,y(x)\}$ with $|y'|\le 1/10$, $|y''|\le 1/10$. Moreover,
$z^{(k)}$ is the unique element of $\mathcal C^{(k)}$ in $\omega$
(we really mean the iterate $\Phi^\ell$ with $\ell=$ integer part
of $\theta k$, but do not want to overload the notations). Notice
that this condition is satisfied for \( k=0,1 \) since  \( W \)
can be chosen to be admissible for \( b \) sufficiently small.

The fact that each critical point has some minimum space
on either side, on which no other critical points lie, and that the
critical points up to order \( k \) must lie on a piece of \(
W^{(u)}(p) \) of finite length (since they are of generation \(
\leq \theta k \)) implies a bound on the possible number of
critical points of order \( k \). This bound will be made explicit
below and will play an important part in the estimates.

\subsubsection{Bound neighbourhoods}\label{sss.quadratic}
For all $n-1\geq k\geq j  \geq 0 $ and $z^{(k)}\in\cal Ck$ we let
$z_j^{(k)}=\Phi^{j+1}(z^{(k)})$ and
\[
B^{(j)}(z^{(k)}_{0}) = \{\xi_{0} :  |\xi_{i}- z_{i}^{(k)}|\leq
e^{-2\alpha i} + 10^{-k} \text{ for all } i\in [0,j]\}.
\]
This is a way of formalising the idea that there is a set of
points which \emph{shadow}, or remain
 \emph{bound} to, the
orbit of $z^{(k)}$ up to time $j$.
The sequence of iterates $0, \dots, k$ is divided into
\emph{free iterates} and \emph{bound iterates}: $j$ is bound if it
belongs to the binding period associated to a return $\nu < j\leq  k$,
i.e.  all the points bound to \( z^{(k)} \) up to time \( k \) are
also bound to another critical point \( \zeta \) between the iterates
\( \nu+1 \) and \( j \). This is explained precisely in
Section \ref{sss.bindingperiods} below. If \( j \) is not a bound
iterate, it is called  a free iterate. By
convention $0$ is a free iterate.

\subsubsection{Hyperbolicity and distortion}
We assume that
the differential map $D\Phi^{k}$ satisfies uniformly hyperbolic
estimates on the bound neighborhood of every $z^{(k)}$:
\begin{equation*}\tag*{\EG{k}}
\|w_{j}(\xi_0)\| :=\|D\Phi^j_{\xi_0}(w_{0})\|\geq e^{\kappa j}.
 \end{equation*}
for all $j\in [0,k]$, any $\xi_0\in B^{(k)}(z_0^{(k)})$, and any
tangent vector $w_{0}$ with slope $\le 1/10$. In particular the
stable and unstable foliations \( \mathcal E^{(j)}, \mathcal
F^{(j)} \) are defined in the whole of \( B^{(j)}(z^{(k)}_{0}) \)
and, as part of the inductive assumptions, the leaves of $\cal Fk$
are admissible curves. Let $f^{(j)}$ be a norm $1$ vector field
tangent to the leaves of \cal Fj respectively, and $f^{(j)}_{i}$
its image under the differential map $D\Phi^{i}$, for $i\ge 0$.
For every free iterate $j < k$ of the critical point
$z^{(k)}\in\cal Ck$, the hyperbolic coordinates \( \cal Hj \) as
well as their images $\mathcal H^{(j)}_{j}$ are $C^{2}$ close to
the standard coordinate system \( \cal H0 \) (i.e. the unstable
leaves are admissible) and satisfy uniform distortion bounds:
there exists a constant $D_0>0$ such that for all points $\xi_{0},
\eta_{0} \in B^{(k)}(z_{0}^{(k)})$ we have
\begin{equation*}\tag*{\BD{k}}
\begin{array}{c}
    \log \frac{\|f^{(j)}_{j} (\xi_{0}) \|}{\|f^{(j)}_{j}(\eta_{0}) \|}
    \leq D_0 \sum_{i < j} \frac{|\xi_{i}-\eta_{i}|}{\ega i}
    \\ \text{and}\\
    \ang (f^{(j)}_{j} (\xi_{0}), f^{(j)}_{j}(\eta_{0}) )
    \leq D_0 \frac{|\xi_{j}-\eta_{j}|}{\ega j}\,.
\end{array}
\end{equation*}
The bounded distortion property says that the orbits of all these
points are in a sense indistinguishable from an analytic point of
view.

\subsubsection{Quadratic tangencies}
As mentioned above, the exponential growth condition guarantees
that hyperbolic coordinates \( \mathcal H^{(k)} \) are defined in
the whole of the bound neighbourhood \( B^{(k)}(z^{(k)}_0)  \).
The critical point \( z^{(k)} \) is then characterized by the
property that the corresponding critical value \( z_{0}^{(k)} \)
is a point of tangency between the image \( \gamma_{0} \) of the
admissible curve \( \gamma \) containing \( z^{(k)} \) and a the
leaf \( \mathcal E^{(k)}(z_{0}^{(k)})\) of the stable foliation of
order \( k \). Again it is possible to show that the curvature of
\( \gamma_{0} \) is much larger than the curvature of the stable
leaves of \( \mathcal E^{(k)} \) and thus this tangency is unique
and quadratic.

\subsubsection{Nested neighborhoods and ancestors}\label{sss.nested}

There exists a sequence $$z^{(k-1)}, \dots, z^{(0)}$$ of
\emph{ancestors} of $ z^{(k)} $ such that for $ i=0, \dots, k-1 $
we have  $ z^{(i)}\in\cal Ci $  and
\[
B^{(i+1)}(z^{(i+1)}_0)\subset B^{(i)}(z^{(i)}_0).
\]
Notice that the term $10^{-k}$ in the definition of the bound neighbourhoods
is much smaller than $e^{-2\alpha i}$, for any $i\le k$, and so it
is negligible from a geometrical point of view. It is introduced
for formal reasons only, to ensure this nested property, see Section
\ref{old}.

\subsubsection{True critical points}
The set of critical points $ \mathcal C $ is obtained as the set
of limit points of any sequence $ \{z^{(k)}\} $ with $
z^{(k)}\in\cal Ck $ and such that $ z^{(k-1)}, \dots, z^{(0)} $
are ancestors of $ z^{(k)} $ for each $ k $.

\subsection{Bounded recurrence}\label{sss.binding}

We assume that the critical sets \( \mathcal C^{(k)} \) are defined
and satisfy the conditions stated above for all \( k\leq n-1
\) and explain how to formulate a bounded recurrence condition
\( (*)_{n-1} \) on the set \( \mathcal C^{(n-1)} \).

\subsubsection{The recurrence condition}\label{recurrencecondition}

Let $0 \leq \nu \leq k$ be a free iterate and $z^{(\nu)}$ be the
corresponding ancestor of $z^{(k)}$. If the $\nu$'th image of
$B^{(\nu)}(z^{(\nu)}_0)$ intersects $\Delta$ we say that $ \nu $
is a \emph{free return} for $z^{(k)}$. Then there is an algorithm,
the \emph{capture argument}, which associates to \( z^{(k)}_{\nu}
\) a particular critical point \( \zeta^{(\nu)}\in\mathcal
C^{(\nu)} \) in \emph{tangential position} to it. We just give a
snapshot of this algorithm at time $\nu$, referring the reader to
\cite{BenCar91}*{\S~6} or \cite{MorVia93}*{\S~9} for the detailed
construction.

\begin{figure}[h]
\begin{center}
 \psfrag{zn}{{\footnotesize $z^{(k)}_{\nu}$}}
 \psfrag{h}{{\footnotesize $\ge e^{-\alpha \nu}$}}
 \psfrag{v}{{\footnotesize $\le b^{c g_j}$}}
 \psfrag{z0}{{\footnotesize $\zeta_{(0)}$}}
 \psfrag{z1}{{\footnotesize $\zeta_{(g_j)}$}}
 \includegraphics[width=4in]{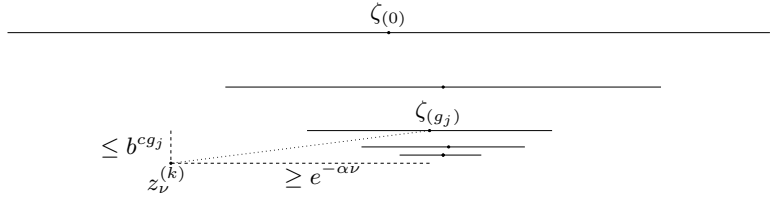}
 \caption{\label{capture} Looking for a binding point}
\end{center}
\end{figure}

As part of the argument, one constructs a whole sequence of
candidates $\zeta_{(g_j)}$ which are critical points sitting on
admissible segments of radii $\rho^{g_j}$ inside $\Phi^{g_j}(W)$
such that the \emph{vertical} distance between \(z^{(k)}_{\nu} \)
and $\zeta_{(g_j)}$ is \(\leq b^{c g_{j}} \) as shown in
Figure~\ref{capture}.
These points are defined for an increasing sequence of $g_i$ which
is not too sparse: $g_{i+1}\le 3 g_i$\,. Then one chooses as the
binding point $\zeta^{(\nu)}=\zeta_{(g_j)}$ where $g_j$ is largest
such that $\zeta_{(g_j)}$ is defined and (the constant
$\theta=C/|\log b|$ was introduced in Section~\ref{sss.constants})
\begin{equation}\label{eq.generationnext}
g_j \le \theta \nu\,.
\end{equation}
This condition will be explained in Section~\ref{tangential}.

Then we define the ``distance'' of \( z^{(k)}_{\nu} \) from the
critical set \( \mathcal C^{(\nu)} \) as the minimum distance
between \( \xi_{\nu} \) and \( \zeta^{(\nu)} \) over all points \(
\xi_{\nu} \) where \( \xi_{0}\in B^{(k)}(z^{(k)}_0) \):
$$
d(z^{(k)}_{\nu}) := |z^{(k)}_{\nu}-\mathcal C^{(\nu)}| :=
\min_{\xi_{0}\in B^{(k)}(z^{(k)}_{0})} |\xi_{\nu} - \zeta^{(\nu)}|
$$
With the notion of distance to the critical set defined above, we can
formulate precisely the bounded recurrence condition
\begin{equation*}\tag*{\( (*)_{k} \)}
  \sum_{\substack{\text{free returns}\\ \nu\leq k}}\!\!\!
 \log d(z_{\nu}^{(k)})^{-1} \leq \alpha k.
  \end{equation*}
Notice that this implies in particular \( d(z^{(k)}_{\nu}) \geq
e^{-\alpha \nu}  \) and even $d(\xi_{\nu})\ge e^{-\alpha\nu}$  for
all \( \xi_{0}\in B^{(k)}(z_{0}^{(k)}) \) and for all $\nu\le k$.
We assume that all critical sets \( \mathcal C^{(k)} \) satisfy
condition \(
 (*)_{k} \) (as well as \EG{k}, \BD{k} and the other conditions
 given above, for all \( k\leq n-1 \) and   prove that this implies
 that conditions \( \EG n \) and
\BD n hold.

\subsubsection{Tangential position}\label{tangential}

A key consequence of the bounded recurrence condition and the
capture argument outlined in Section~\ref{recurrencecondition} is
that the binding point $\zeta^{(k)}$ and $\xi_k$ are in
\emph{tangential position} for all $ \xi_0\in B^{(k)}(z^{(k}_0)$:
there exists an admissible curve \(\gamma \) which is tangent to
the vector $w_{k}(\xi_0)$ at $\xi_{k}$ and tangent to the unstable
manifold at $\zeta^{(k)}$. In particular, the critical point
\(\zeta^{(k)}\) chosen via the capture argument has essentially
the same vertical coordinate as any of these \( \xi_{k} \)\,,
including the critical iterate \( z^{(k)}_{k} \)\,.

Indeed, the bounded recurrence condition \( (*)_{k} \) implies
that the horizontal distance from $z_k$ to the binding point is
$\ge e^{-\alpha k}$. Hence, to ensure tangential position we have
the choice of any $\zeta_{(g_j)}$ with
$$
e^{-\alpha k} \gg b^{c g_j} \quad\text{or equivalently} \quad
g_j\ge\frac{\const}{|\log b|}\, k.
$$
This shows, in other words, that it is sufficient to consider
critical points of generations \( g\le\const k/|\log b| = \theta k
\) to guarantee the existence of one in tangential position.
Indeed, this is how the expression $\theta=\theta(b)$ and
condition \eqref{eq.generationnext} come about. As we shall see
below, this also guarantees that the number of critical points of
a given order are not too many to destroy the parameter exclusion
estimates.

The reason being in tangential position is so crucial is that it
allows for estimates at returns which are very much the same as in
the one-dimensional situation. In particular, the ``loss of
expansion'' is roughly proportional to the distance to the binding
critical point. See \cite{BenCar91}*{\S~7}, \cite{MorVia93}*{\S~9}
and Section~\ref{localgeometry} below.

\begin{remark}
In \cites{BenCar91, MorVia93} the critical set is constructed in
such a way that the tangential position property at free returns
is satisfied for the critical points themselves. One main
contribution in \cite{BenYou93} was to show that the capture
argument works for essentially any other point in \( W^{(u)}(p) \)
as well and this implies the existence of a hyperbolic
Sinai-Ruelle-Bowen measure. Further results such as exponential
decay of correlations \cite{BenYou00} and other hyperbolicity and
topological properties \cite{WanYou01} ultimately rely on this
fact. Moreover, \cite{BenVia01} went one step further and proved
that for Lebesgue almost all points \emph{in the basin of
attraction} returns are \emph{eventually} tangential. This is
crucial in their proof that the basin has ``no holes'': the time
average of Lebesgue almost every point (not just a positive
measure subset) coincides with the Sinai-Ruelle-Bowen measure.
\end{remark}

\subsection{Hyperbolicity and distortion at time \protect\(n \protect\).}
\label{ss.main}

We outline the proof of the first part of Theorem \ref{positive
exponents} where the bounded recurrence condition on the critical
set \( \mathcal C^{(n-1)} \) is shown to imply some hyperbolicity
and distortion estimates in a neighbourhood of each point of \(
\mathcal C^{(n-1)} \) up to time \( n \). The situation we have to
worry about is  when the critical point has a return at  time  \(
n-1 \), otherwise the calculations are relatively straightforward.
In the case of a return however, as we mentioned above, vectors
get rotated and end up in almost vertical directions which are
then violently contracted for  many iterations, giving rise to a
possibly unbounded loss of expansion accumulated up to time \( n
\). The idea therefore is to use condition \( (*)_{n-1} \) to
control the effect of these returns. We  assume that \( n-1 \) is
a free return as in the previous section and let \(
\zeta=\zeta^{(n-1)} \) denote the corresponding associated
critical point.

\subsubsection{Binding periods}\label{sss.bindingperiods}
\label{length}

Our inductive assumptions imply that hyperbolic coordinates \(
\mathcal H^{(k)} \) are defined in the neighbourhoods \(
B^{(k)}(\zeta_{0}) \) for all \( k\leq n-1 \). Notice that these
bound neighbourhoods shrink as \( k \) increases, but start off
relatively large for small values of \( k \). Therefore there must
be some values of \( k \) for which the image \( \gamma_{0} =
\Phi(\gamma) \in  B^{(k)}(\zeta_{0}) \). We denote by \( p \) the
\emph{largest} such \( k \). In principle we do not know that \( p
< n-1 \) but it is not difficult to prove that in fact that $p
\sim |\log d(z_{n-1}^{(n-1)})|$. In particular, under condition
$(*)_{n-1}$ we get $p\sim \alpha n \ll n$ (fixing $\alpha$ small).
Thus the point \( \xi_{n} \) will \emph{shadow} \( \zeta_{0} \)
for exactly \( p \) iterations. We say that \( p \) is the length
of the \emph{binding period} associated to the return of \( z_{0}
\) to \( \Delta \) at time \( n-1 \).

\subsubsection{Local geometry}\label{localgeometry}

We now want to analyse carefully
the geometry of \( \gamma_{0} \) and
\( w_{n}(\xi_0) \) with respect to the hyperbolic coordinates
\( \mathcal H^{(p)} \). The information we have is
that \( \gamma_{0} \)
is tangent at \( \zeta_{0} \) to a stable leaf  \( \mathcal
E^{(n-1)}(\zeta_{0}) \) and that the curve \( \gamma_{0} \) is
quadratic with respect to the coordinate system \( \mathcal H^{(n-1)}
\). Now suppose for the moment that \( p=n-1 \).  Then the
quadratic nature of \( \gamma_{0} \) with respect to \( \mathcal
H^{(n-1)} \) implies that the \emph{slope} of \( w_{n}(\xi_0) \)
\emph{in these coordinates}
is related to the distance between \( \xi_{n-1} \) and \( \zeta \)
and more specifically the ``horizontal'' component of \( w_{n}(\xi_0)
\),  that is, the component in the direction of \( f^{(n-1)}(\xi_{n}) \)
is proportional to \( d(z_{n-1}^{(n-1)}) \).

These estimates do not apply immediately to hyperbolic
coordinates for arbitrary \( p\leq n-1 \), for example they may not
apply to the standard coordinates \( \mathcal H^{(0)} \) as the
\( w_{n}(\xi_0) \) may actually be completely vertical in these
coordinates and therefore have no horizontal component. Nevertheless
it follows from \eqref{anglek} that the angle between leaves
associated to \( \mathcal H^{(n-1)} \) and
\( \mathcal H^{(p)} \) for \( p \leq n-1\) is of order \( b^{p} \).
Moreover \( p\sim \log  d(z_{n-1}^{(n-1)})^{-1} \ll \log
d(z_{n-1}^{(n-1)})/\log b \)
and therefore  \( b^{p}\ll  d(z_{n-1}^{(n-1)}) \) and therefore
the length of the horizontal component in the coordinates \( \mathcal
H^{(p)} \) is essentially the same in \( \mathcal H^{(n-1)} \).

\subsubsection{Recovering hyperbolicity}

The fact that the horizontal component of \( w_{n}(\xi_0) \) (in
hyperbolic coordinates) is proportional to \( d(z_{n-1}^{(n-1)})
\) is a two-dimensional analogue of the simple fact that in the
one-dimensional case, the loss of derivative incurred after a
return to \( \Delta \) is proportional to the distance to the
critical point. Thus, even though the vector \( w_{n}(\xi_0) \)
may be very close to vertical (in fact it may be vertical in the
standard coordinates) and therefore suffer strong contraction for
arbitrarily many iterates, we do not need to worry about the
contraction because we know that it has a component of strictly
positive length proportional to \( d(z_{n-1}^{(n-1)}) \) and thus
of the order of \( e^{-\alpha n} \) in the ``horizontal''
direction and this component is being expanded, providing us with
a lower bound for the real size of the vector. Using the inductive
assumptions we can show that an average exponential rate of growth
is recovered by the end of the binding period:
\begin{equation}\label{eq.growthinbindingperiods}
\|D\Phi^{p+1}(\xi_{n-1}) w_{n-1} \|
 \ge d(\xi_\nu)^{-\kappa_1}
 \ge e^{\kappa_2 p} \gg 1,
\end{equation}
for all $\xi\in B^{(n+p)}(z^{(k)}_{0})$, where the constants
$\kappa_1$, $\kappa_2>0$ depend only on $\kappa$. In fact the
strong contraction is useful at this point because it implies that
the ``vertical'' component, i.e. the component in the direction of
\( e^{(n-1)}(\xi_{n}) \) is shrinking very fast and this implies
that the slope of the vector is decreasing very fast and that it
returns to an almost horizontal position very quickly.

\subsubsection{Bounded distortion}

Using the geometrical structure and estimates above one also
proves that the bounded distortion property \BD n holds. This is a
technical calculation and we refer to \cites{BenCar91, MorVia93}
or \cite{LuzVia} for the proof in much the same formal setting as
that given here.

\subsection{New critical points}\label{ss.new critical points}

We give two algorithms for generating the new critical set \(
\mathcal C^{(n)} \). Both of them depend on the fact that
since condition \EG n is satisfied by all points of \( \mathcal
C^{(n-1)} \) it follows in particular that the hyperbolic coordinates
\( \mathcal H^{(n)} \) are also defined in neighbourhoods of these
points.

\subsubsection{Refining the set of critical points of order \( n-1 \)}
\label{old}

Since \( \mathcal H^{(n)} \) does not generally coincide with \(
\mathcal H^{(n-1)} \), the  critical points \( z^{(n-1)} \) are no
longer tangent to the new stable foliations \( \mathcal E^{(n)}
\). Instead, these foliations define new points of tangencies with
the new stable leaves \( \mathcal E^{(n)} \) close to the old
ones. By definition these belong to the new set \( \mathcal
C^{(n)} \) of critical points of order \( n \). By the estimates
on the convergence of hyperbolic coordinates, see e.g.
\eqref{anglek}, the ``distance'' between the leaves of \( \mathcal
E^{(n-1)} \) and the leaves of \( \mathcal E^{(n)} \) is of the
order \( b^{n-1} \) and therefore the distance between the new
points of tangencies, i.e. the new critical points, and the old
ones will also be of the order of \( b^{n-1} \), which is
extremely small. It is then easy to see that the distance between
the iterates \( z_{j}^{(n-1)} \) and \( z_{j}^{(n)} \) will
continue to be essentially negligible for all \( j\leq n \). In
particular the nested property of bound neighbourhoods is
satisfied, and the new point \( z^{(n)} \) inherits all the
properties of its ancestor \( z^{(n-1)} \) as far as bounded
recurrence, exponential growth, and bounded distortion are
concerned.

\subsubsection{Adding really new critical points}\label{new}

Notice that there may be other admissible pieces of the unstable
manifold \( W^u(P) \) which are too small or not on the right
section of \( W^u(P) \) to admit critical points of order \( n-1
\), recall property \ref{sss.admissible}, but can in principle
admit critical points of order \( n \). We add these points to the
new critical set \( \mathcal C^{(n)} \) as long as they are close
enough to \( \mathcal C^{(n-1)} \) so that in particular the
nested property of bound neighbourhoods is satisfied. This
completes the definition of \( \mathcal C^{(n)} \) and the sketch
of the proof of Theorem \ref{positive exponents}.

\subsection{The cardinality of the critical set}\label{cardinality}

Before going on to discuss the parameter dependence of the objects
defined above, we make a couple of important remarks regarding the
definition of the set \( \mathcal C^{(n)} \).

\subsubsection{Why we need many critical points}

We recall that the overall objective of our discussion is to prove
the existence of many parameters for which some (non-uniform)
hyperbolicity conditions are satisfied. As a first step in this
direction, it is useful to start with the relatively modest
objective of showing that the  unstable manifold \( W^{(u)}(p) \)
is not contained in the basin of attraction of an attracting
periodic orbit, a necessary, though not sufficient, condition for
the hyperbolicity conditions to hold. To prove this it is enough
to show that almost all points \( z\in W^u(P) \) satisfy the
exponential growth condition \EG n for all \( n\geq 1 \). The
proof of this fact requires controlling returns to \( \Delta \)
and the argument presented here relies on achieving this control
by identifying a set of \emph{critical points} as explained above,
with the crucial property that a critical point in tangential
position can always be found at every free return as long as the
bounded recurrence condition is satisfied. This critical point can
then be used  to implement the shadowing (binding) argument to
show that the exponential growth condition can be maintained
through the passage in \( \Delta \). Since returns can occur at
various ``heights'', tangential position can only be guaranteed if
there are \emph{sufficiently many} critical points.

\subsubsection{Why we need not-too-many critical points}

A choice of critical set which contains many points becomes
problematic in view of our strategy of defining a good parameter
in terms of some recurrence conditions on such points. The more
critical points there are the greater the likelihood that at least
one of them will fail to satisfy such condition and will lead to
having to exclude a particular parameter value. Therefore, it is
crucial to ensure that there are \emph{relatively few} critical
points such that by imposing the recurrence condition on their
orbits one controls the whole dynamics, in the sense that one si
able to prove hyperbolicity. Ultimately, at least at the present
stage of the theory, this requires a strong (smallness)
restriction on the perturbation size $b$.

\subsubsection{A reasonable compromise} \label{sss.afford1}

The main restriction on the number of critical points of a given
order \( k \) comes from the requirement that they are of
generation \( g \le \theta k \) and that they have some space
around them where there is no other critical point, see Section
\ref{sss.admissible}. These properties immediately imply the
following crucial bound on the total number of critical points of
order $k$:
\begin{equation}\label{eq.numberofcriticalpoints}
\#\mathcal C^{(k)}
 \le \frac{|\Phi^{\theta k}(W)|}{2\rho^{\theta k}}
 \le (5/\rho)^{\theta k}.
\end{equation}
The constant $5$ is an upper bound for the norm of the derivative.

We shall see in the parameter exclusion argument that this bound
is good enough to ensure that not too many parameters get
excluded. On the other hand, the reason we can afford to use only
critical points with the above properties is related to the
features of the constructions in
Sections~\ref{recurrencecondition} and \ref{tangential}: as we
have seen, a binding critical in tangential position can always be
found among the critical points of generation \( g \leq \theta k
\) and lying on admissible unstable segments of radius
$\rho^g\ge\rho^{\theta k}$.

\section{Positive measure in parameter space}\label{s.positivemeasure}

Next we explain why the set of parameter values $a$ for which the
previous construction works has positive Lebesgue measure. It is
assumed that \( b \) is sufficiently small and that
$a$ varies in an interval $\Omega$ close to $a=2$ and
not too small.

\begin{theorem} \label{positive measure}
There exists a set
 $\Gamma^{*}\subset\Omega $ such that:
\begin{enumerate}
\item the Lebesgue measure $ |\Gamma^{*}| >0 $;
\item for all $ a\in\Gamma^{*} $
 a critical set $ \mathcal C $ is  defined
and satisfies $ (*)_{n} $ for all $ n\geq 0 $.
\end{enumerate}
\end{theorem}

The precise condition for the choice of the interval \( \Omega \)
is in terms of the limiting
one-dimensional map $\phi_a(x)=1-a x^2$. Firstly, the iterates
$c_n$ of the critical point remain outside the critical region
$\Delta$ for the first $N$ iterates, for some large $N$. Secondly,
$c_N(a)=\phi_a^{N+1}(c)$ describes an interval of length
$>\delta/10$ in a monotone fashion when $a$ varies in $\Omega$.
This last requirement ensures that $N$ is an \emph{escape
situation} (this notion will be recalled in a while). By simple
perturbation, these properties extend to the two-dimensional
H\'enon-like map $\Phi_a$ if $b$ is sufficiently small.

The proof relies on the construction of a nested sequence of sets
$ \Gamma^{(n)} $ such that each parameter value in $ \Gamma^{(n)}
$ has a critical set $ \cal Cn $ satisfying $ (*)_{n-1} $.
 The set $ \Gamma^{*} $ is then just the intersection of all
$ \Gamma^{(n)} $.
 The main estimate concerns the probability of exclusions at each
time $n$, that is, the Lebesgue measure of
 $ \Gamma^{(n-1)}\setminus\Gamma^{(n)} $.
 We begin here with a sketch of the construction of the sets $
\Gamma^{(n)} $ in the one-dimensional case and discuss the main
issues with the generalizations of the construction to the
two-dimensional setting.

\subsection{The one-dimensional case}\label{ss.one-dimensional}

Given the critical point $z=c=0$ and an integer $k\geq 0$ we
define the map
\begin{equation}\label{parameter maps}
z_{k}:\Omega \to Q,
 \qquad z_{k}(a) = \phi_{a}^{k+1}(z(a)),
\end{equation}
from parameter space to phase space associating the $k$:th iterate
of the critical value $z_{0}(a)=\phi_a(z)$ to each parameter value
$a\in\Omega$.
Whenever $z_{k}(\Omega)$ intersects the critical neighborhood
$\Delta$ we subdivide it into subintervals by pulling back a
certain partition $\mathcal I$ of $\Delta$. Roughly, the partition
consists of the intervals bounded by the sequence $\pm e^{-r}$ for
$r\ge |\log\delta|$ (for distortion reasons these intervals must
be subdivided a bit further). Then we exclude those parameter
subintervals for which condition $(*)$ does not hold at time $k$.

We obtain in this way a sequence of good parameter sets
$\Gamma^{(n-1)}\subset\dots\subset\Gamma^{(0)}=\Omega$ and
corresponding partitions $\cal P{n-1}, \dots, \cal P0$ such that
all parameters in any given $\gamma\in\cal Pk$ have essentially
indistinguishable itineraries (in particular as far as the
critical recurrence is concerned) and essentially equivalent
derivative estimates  up to time $k+1$ (in particular $z_{k}$
restricted to elements of $\cal Pk$ is a diffeomorphism onto its
image).

At each step we refine $\cal P{k}$ to a partition $\hcal P{k}$ of
$\Gamma^{(k)}$ by pulling back the intersection of elements of
$\cal P{k}$ under the map $z_{k+1}$ with $\mathcal I$. We then
exclude those elements of $\hcal P{k}$ for which the recurrence
condition fails, and define $\Gamma^{(k+1)}$ as the union of the
remaining elements and $\cal P{k+1}$ as the restriction of $\hcal
P{k}$ to $\Gamma^{(k+1)}$. A large deviations type of argument
shows that the measure of the excluded set decreases exponentially
fast with $k$:
$$
|\Gamma^{(k)}\setminus\Gamma^{(k+1)}| \le e^{-\tau_0 k} |\Omega|
\quad\text{for all } k\ge N,
$$
where $\tau_0>0$ is independent of $N$. Taking $N$ large enough
(no exclusions are needed inside $\Omega$ before time $N$), this
gives that a positive measure set remains after all exclusions.
We do not give the details here as this is a special case  of the
argument in the two-dimensional context, which will be discussed in
some detail below.

\subsection{Two-dimensional issues}
We mention here the key differences between the one-dimensional
and two-dimensional situations and the main difficulties in
generalizing the scheme sketched above to the two-dimensional
case.

\subsubsection{Many critical points}
The most obvious difference is that in two dimensions there is a
large number of critical points of order \( n \) at each stage $n$
and \emph{all} these critical points must satisfy $(*)_{n}$. Thus
many more parameter exclusions are necessary. However we have seen
in (\ref{eq.numberofcriticalpoints}) that the cardinality of
${\mathcal C}^{(n)}$ grows at most exponentially fast with $n$,
with exponential growth rate which can be made arbitrarily small
by reducing $b$. This is crucial to guarantee that the total
proportion of parameters excluded at time $n$ continues to be
exponentially small in $n$: the measure of parameters excluded by
imposing the recurrence condition on each individual critical
point decreases exponentially fast with $n$, with decay rate
$\tau_0$ which is essentially the same as in dimension one, and so
is independent of $b$. Section~\ref{ss.mainmeasureestimate} makes
these explanations more quantitative.

\subsubsection{Interaction between different critical points}

A second important issue is that the one dimensional argument
relies on keeping track of combinatorial and analytic data related
to the history of the critical orbit for various parameter values.
Here we can do the basically the same, but each one of the
critical orbits has its own associated data, since the dynamical
history and pattern of recurrence to the critical neighbourhood
vary from one critical orbit to the other. For this reason, it
will be convenient to introduce an extended parameter space, with
separate combinatorial structures (partitions, itineraries)
relating to each critical point. While we try to think of these
structures as being essentially independent, this is not entirely
accurate because different critical orbits do interact with each
other. Namely, a critical point $z$ may require a different one
$\zeta$ as the binding point associated to some free return. Then,
if a parameter is deleted because $\zeta$ fails to satisfy
condition $(*)_{n}$ for that parameter value, this deletion must
be somehow registered in the combinatorial structure of the other
critical points $z$. Section~\ref{ss.interaction} explains how
this is handled.

\begin{remark}
Neither of these two points is really related to dimensionality:
multiplicity of critical points and the difficulties connected to
interactions between their orbits occur already for
\emph{multi}-modal maps in dimension $1$. In fact,
\cite{PacRovVia98} treated those difficulties in the extreme case
of \emph{infinite}-modal maps of the interval, that is with
infinitely many critical points, using this strategy of defining
different but not-quite-independent combinatorics and exclusion
rules for each critical point that we just outlined and will be
detailing a bit more in a while.
\end{remark}

\subsubsection{Continuation of critical points}
Another fundamental difficulty, this time intrinsically
two-dimensional, is the problem of talking about a given critical
point for different parameter values, as was implicitly assumed in
the discussion of the previous two points. It is not immediately
obvious how to do this because critical points are defined
dynamically: the definition requires certain hyperbolicity
properties to be satisfied and the precise location of the point
depends on the geometrical and dynamical features of the map for a
specific parameter value, which are very unstable under parameter
changes. We shall use the fact that critical points of finite
order do admit a \emph{critical continuation} to a neighborhood in
parameter space: the condition of quadratic tangency that defines
such points has a unique smooth solution on that neighborhood. As
a matter of fact, we make it here an additional requirement on a
tangency of order $k$, for it to be in the critical set $\cal
Ck_a$, that it should have a suitable continuation in parameter
space.

To appreciate the situation better, suppose for example that a
critical point of order \( k \) admits a critical continuation to
a parameter interval \( \omega \). Suppose however that there
exists a subinterval \( \tilde\omega\subset\omega \) such that \(
\omega\setminus\tilde\omega \) has two connected components and
such that the required bounded recurrence condition \( (*)_{k} \)
fails to be satisfied by the critical point for \(
a\in\tilde\omega \). Then the critical point cannot be refined to
an approximation of order \( k+1 \) for critical points in \(
\tilde\omega \) although it can in the two components of \( a\in
\omega\setminus\tilde\omega \). We need to address the questions
of whether these refinements can still be thought of as
continuations of each other, i.e. whether we can still talk about
a single critical point with a critical continuation on the
(disconnected) set \( \omega\setminus\tilde\omega \) or whether we
should think of having two independent critical points defined in
the two distinct parameter intervals. See
Sections~\ref{ss.mainmeasureestimate} and \ref{ss.extended} for
the details of how these issues are resolved and how we manage to
relate critical points existing for different parameter values.

\subsection{Overview of the argument}

All of these issues will be dealt with formally by defining an
\emph{extended parameter space} where each critical point (of
finite order) comes with its own interval of parameters on which
it admits a continuation as a critical point, and with its own
combinatorial and analytical data. In the remaining part of this
section we describe the structure of this extended parameter space
at each iterate $k$, and outline the main calculation that proves
that $\Gamma^{*}=\cap_k \Gamma^{(k)}$ has positive Lebesgue
measure.

\begin{figure}[h]
\begin{center}
 \psfrag{a}{$a$}
 \psfrag{G}{$\Omega_{[z]}$}
 \psfrag{O}{$\Omega$}
\includegraphics[width=4in]{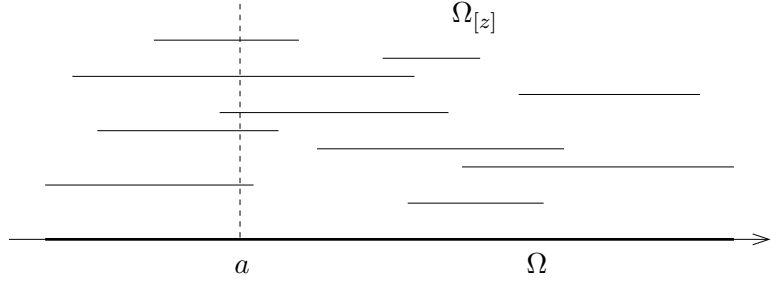}
\caption{\label{parameter space} Extended parameter space}
\end{center}
\end{figure}

\subsubsection{The extended parameter space}
\label{ss.mainmeasureestimate}

The parameter space at time $k$ consists of a disjoint union of
copies $\Omega_{[z]}$ of (not necessarily disjoint) subintervals
of $\Omega$, as described in Figure~\ref{parameter space}. Each of
them comes with a critical point $z^{(k)}(a)$ of order $k$ defined
on some subset of \( \Gamma^{(k)}_{[z]} \) of $\Omega_{[z]}$. The
symbol \( [z] \) parametrizes the set of these segments
$\Omega_{[z]}$ and may be thought of as an ``equivalence class''
of critical points in the sense that there exists some \( \nu\leq
k \) and a critical point $z^{(\nu)}(a)$ which admits a
continuation as a tangency of order $\nu$ over the entire
$\Omega_{[z]}$ and which is ancestor to $z^{(k)}(a)$ whenever the
latter is defined. For this reason, it makes sense to think of
$z^{(k)}(a)$ as ``the same critical point''  for different
parameter values in its domain. In addition, $z^{(\nu)}$ has an
{\em escape situation\/} at time $\nu$~: the image of
$$\Omega_{[z]}\ni a \mapsto z_\nu^{(\nu)}(a)$$ is an
admissible curve of length $\ge\delta/10$. This exactly
corresponds to the requirement, in the one-dimensional setting,
that the initial parameter interval should not be too small.

The set \( \Gamma^{(k)}_{[z]} \) is a finite union of subintervals
of \( \Omega_{[z]} \) and is a subset of parameters in \(
\Omega_{[z]} \) for which the corresponding critical point \(
z^{(k)}(a) \) satisfies the recurrence conditions up to time $k$.
It also comes with a combinatorial structure in the form of a
finite partition \( \mathcal P^{(k)}_{[z]} \) into subintervals
defined in such a way that all critical points \( z^{(k)}(a) \)
with \( a \) belonging to any one element of this partition have
the same history, that is, essentially the same analytic,
hyperbolicity, distortion, and recurrence estimates, up to time \(
k-1 \).

\subsubsection{The parameter exclusion argument}
\label{sss.parexc} For each $[z]$ we exclude a set of parameters
$E^{(k)}_{[z]}\subset \Omega_{[z]}$ to enforce condition
$(*)_{k}$. These individual exclusions are estimated in much the
same way as in dimension $1$, the details will be given in the
following sections. For the moment we just mention that we begin
by defining a refined partition \( \hcal P k_{[z]} \succ \mathcal
P^{(k)}_{[z]} \) of the parameter set \( \Gamma_{[z]}^{(k)} \),
depending on the position of the critical points \( z^{(k)}_{k}(a)
\) for each parameter \( a\in \Gamma_{[z]}^{(k)} \). We then
decide which parameters to exclude on the basis of this additional
combinatorial information. We always exclude whole elements of
this refined partition, and not just individual parameters, even
if this may mean excluding somewhat more parameters than is
actually necessary. This is important because it ensures that the
remaining set
$$
\Gamma^{(k+1)}_{[z]} = \Gamma^{(k)}_{[z]}\setminus E^{(k)}_{[z]}
$$
of parameters which are good for \( [z] \) up to time \( k \)
inherits a combinatorial structure, the family of atoms of the
refined partition which have not been excluded, and these are
relatively large intervals. Indeed, our exclusion estimates depend
crucially on lower bounds on the size of parameter intervals
(small intervals might even be completely deleted at one given
return!), and removing individual parameter values could lead to
the formation of such small connected components in parameter
space. In Section~\ref{ss.conclusion} we get
\begin{equation}\label{meas2}
|E^{(k)}_{[z]}| \le e^{-\tau_0 k} |\Omega_{[z]}|
\end{equation}
with $\tau_0>0$ independent of $b$ and $N$. By definition, the new
set of good parameters is
\begin{equation}\label{meas3}
\Gamma^{(k+1)}
 = \Gamma^{(k)} \setminus \bigcup_{[z]}E^{(k)}_{[z]}.
\end{equation}
This means \emph{we only consider a parameter value good at any
given time $k+1$ if it is good for all critical points up to that
time}. Notice also that there are no partitions associated to \(
\Gamma^{k} \) or \( \Gamma^{(k+1)} \), these are just ``raw'' sets
of parameter values.

To estimate the total size of exclusions, we  remark that if a
parameter $a$ belongs to intervals $\Omega_{[z]}$ and
$\Omega_{[w]}$ then, by definition of critical points, the
corresponding ancestor points $z^{(\nu)}(a)$ and $w^{(\mu)}(a)$
are $\ge 2\rho^{\theta k}$ away from each other in the intrinsic
metric of the unstable manifold \( W^{(u)}(p) \). Together with
the fact that they must be contained in a compact part of \(
W^{(u)}(p) \) of length \( \leq 2 \times 5^{\theta k} \) (since we
started with a leaf of length \( \leq 2 \) and iterated this for
at most \( \theta k \) iterates with a maximum expansion of a
factor 5 at each iteration),  the same calculation as for
\eqref{eq.numberofcriticalpoints} gives us the following bound on
the size of this family of intervals (see
Section~\ref{sss.newcritical}):
\begin{equation}\label{meas1}
 \#\{[z]\in\cal Ck: a\in \Omega_{[z]} \}
 \le (5/\rho)^{\theta k} \quad\text{ for any } a\in\Omega.
\end{equation}
So the total exclusions at this iterate are
\begin{equation}\label{meas4}
\big|\bigcup_{[z]}E^{(k)}_{[z]}\big|
 \le \sum_{[z]} \big|E^{(k)}_{[z]}\big|
 \le e^{-\tau_0 k} \sum_{[z]} |\Omega_{[z]}|
 \le e^{-\tau_0 k} (5/\rho)^{\theta k} |\Omega|,
\end{equation}
by (\ref{meas2}) and (\ref{meas1}). Assuming $b$ is small, the
term on the right is $\le e^{-(\tau_0/2) k} |\Omega|$.
\begin{equation}\label{meas5}
|\Omega\setminus\Gamma^{*}|
 = \big|\bigcup_{k=N}^{\infty} \bigcup_{[z]} E^{(k)}_{[z]}\big|
 \le \sum_{k=N}^{\infty} e^{-(\tau_0/2)k} |\Omega|<|\Omega|,
\end{equation}
where $\Gamma^{*}$ is the intersection of all $\Gamma^{(k)}$. The
last inequality assumes $N$ was chosen large enough, and implies
that $|\Gamma^{*}|>0$.

\subsubsection{Interaction between different critical orbits}
\label{ss.interaction}

Observe that each individual parameter interval \(
\Gamma^{(k+1)}_{[z]} \) typically contains some parameter values
which are not in \( \Gamma^{(k+1)} \): at each stage there may
exist (globally) bad parameters which, nevertheless, are good for
\emph{some} of the critical points, at least up to that stage.
This is inevitable, given that we always exclude entire partition
intervals, as explained before, and that different critical points
have different partitions. However, a little bit of thought shows
that this is also most natural to happen.

To explain why, let us consider any parameter value $\bar{a}$ for
which there is a homoclinic point $z$ associated to the fixed
point $P$ (these parameters form a zero measure set, we mention
this situation because it sheds light into the general case). The
forward orbit of $z$ converges to $P$ and, thus, never goes to the
critical region. The recurrence condition is automatically
satisfied, and hyperbolicity features on the homoclinic orbit
follow simply from Proposition~\ref{2dimunifexp}: there is no need
for the binding argument, etc. The point $z$ is a true critical
point (point of tangency between true stable and unstable
manifolds) and from its point of view the parameter $\bar{a}$ is
perfectly good, \emph{notwithstanding the fact that $\bar{a}$ may
be a bad parameter for some other critical point (in which case it
is excluded from $\Gamma^{*}$) and the map may even exhibit
periodic attractors:} this one critical point never becomes aware
of it!

Having said this, different critical orbits do interact with each
other in general. In terms of our inductive construction this
interaction materializes when a critical point $[w]=w^{(k)}$ is
used as the binding point associated to some free return $k$ of a
different critical point $[z]=z^{(k)}$  (this does not occur in
the special situation discussed before): parameters that have been
excluded because $[w]$ failed to satisfy condition $(*)_j$ at some
iterate $j\le k$ must be excluded from the parameter space of
$[z]$ as well. We do indeed exclude an additional set of partition
elements \( \gamma \) of \( \hcal Pk_{[z]} \), but only those
which have already been \emph{completely eliminated} due to
parameter exclusions associated to other critical points, i.e.
such that \( \gamma\cap \Gamma^{(k)} = \emptyset \).  This means
we are really excluding a somewhat larger set
$$
\hat E_{[z]}^{(k)} \supset E_{[z]}^{(k)}
$$
from the parameter interval of $[z]$ at time $k$ and defining
\[
\Gamma^{(k+1)}_{[z]}=\Gamma^{(k)}_{[z]}\setminus \hat E_{[z]}^{(k)}
\subset \Gamma^{(k)}_{[z]}\setminus  E_{[z]}^{(k)}.
\]
An easy, yet
important observation is that these exclusions {\em do not\/}
affect the calculation made before: by definition, any parameter
in the difference
belongs to $E_{[w]}^{(j)}$ for some $[w]$ and some $j\le k$, hence
\begin{equation}\label{justthesame}
\Gamma^{(k)} \setminus \bigcup_{[z]}\hat E^{(k)}_{[z]} =
\Gamma^{(k)} \setminus \bigcup_{[z]}     E^{(k)}_{[z]}.
\end{equation}
In other words, (\ref{meas3}) and (\ref{meas5}) are not changed at
all!

The success of this strategy is based also on the important
observation that \emph{we do not need to exclude an element \(
\gamma \in \mathcal P^{(k)}_{[z]} \) as long as at least one
parameter \( a\in\gamma \) belongs to \( \Gamma^{(k)} \).} This is
explained in more detail in Section
\ref{sss.existenceofbindingpoint} and is essentially due to the
fact that as long as there is even a single parameter \( \tilde a
\in \gamma \cap \Gamma^{(k)} \) then the capture argument works
and there is a binding critical point \( \zeta(\tilde a) \) in
tangential position if \( k \) is a free return for \( \gamma \).

\section{The combinatorial structure}
\label{s.definition}

We are now going to detail the construction outlined in the previous
section.
We begin by giving explicit definitions of the extended
parameter space and the set of good parameters for small values of
$k$.

\subsection{First step of the induction}\label{ss.firststep}

The hyperbolic fixed point $P$ has a continuation $P(a)$ for all
$a\in\Omega$ and we can also consider a continuation $W(a)$ of the
compact interval $W\subset W^u(P)$ introduced in
Section~\ref{ss.properties}. Note that $P(a)$ and $W(a)$ depend
smoothly on the parameter.
Given $1\leq k\leq N$, we have contracting directions $e^{(k)}$ of order
$k$ defined at each point of $W_{0}(a)=\Phi_a(W(a)\cap\Delta)$ and
there exists a unique
point $z^{(k)}(a)\in W(a)$ such that the $k$'th contracting
direction is tangent to $W_{0}(a)$ at the critical value
$z^{(k)}_{0}(a) =\Phi_a(z^{(k)}(a))$.
For $1\leq k \leq N $ we let $\Gamma^{(k)}=\Omega$ and for each
$a\in\Gamma^{(k)}$ we let the critical set $\cal Ck_{a}$ consist
exactly of this critical point $z^{(k)}(a)$. This defines the
critical set $ \cal Ck_{a} $ for all $a\in\Omega$. The extended
parameter space reduces to the single interval
$\Omega_{[z]}=\Omega$.

\subsection{Properties of parametrized critical points}\label{ss.extended}

We now fix $n \ge N$ and suppose inductively that for each $1\leq
k \leq n$ we have already constructed

\begin{itemize}
\item a family of intervals $\{\Omega_{[z]}: [z]\}$ each one
with associated critical point $[z]=z^{(k)}$ defined on a set of
good parameters $\Gamma_{[z]}^{(k)}\subset\Omega_{[z]}$: these
critical points satisfy $(*)_{k-1}$ and \EG k for all $a\in
\Gamma_{[z]}^{(k)}$;
\item and a set $\Gamma^{(k)}$ of parameters good for all critical
points: $\Gamma^{(k)}\cap\Omega_{[z]}$ is contained in
$\Gamma_{[z]}^{(k)}$ for all $[z]$.
\end{itemize}
From now on \( \cal Ck  \) will represent the set of $[z]$
parametrizing the family of intervals above, which we think of as
the set of all critical points or order $k$. To avoid confusion
with the notation below, notice that points \( \cal Ck  \) are
well defined in virtue of their satisfying condition \EG k but are
only assumed to satisfy \( (*)_{k-1} \) (not \( (*)_{k} \)). The
extended parameter space is the disjoint union of intervals:
\[
\Omega_{*}^{(k)}
 = \coprod_{[z]\in\cal Ck}\negfour\Omega_{[z]}.
 \]

These objects have the following additional properties:

\subsubsection{Globally defined ancestor}\label{sss.globally}

There exists a critical point $z^{(\nu)}(a)$ of some order $\nu\le
k$ defined on the whole $\Omega_{[z]}$ which is an ancestor to
$z^{(k)}(a)$ whenever the latter is defined. In addition, $\nu$ is
an escape situation for $z^{(\nu)}$ so that the image of
$\Omega_{[z]}$ under $a\mapsto z^{(\nu)}_{\nu}(a)$ is an
admissible curve with length $\ge \delta/10$.

\subsubsection{Location and uniqueness}\label{sss.location}

Every $z^{(\nu)}(a)$ is the midpoint of an admissible curve
$\omega_a$ of radius $\rho^{\theta \nu}$ inside
$\Phi_a^{\theta\nu}(W(a))$, for $a\in\Omega_{[z]}$. The critical
value $z_0^{(\nu)}(a)=\Phi_a(z^{(\nu)}(a))$ is a point of
quadratic tangency between $\Phi_a(\omega_a)$ and the stable
foliation of order $\nu$ in the bound neighborhood of
$z_0^{(\nu)}(a)$. Moreover, $z^{(\nu)}(a)$ is the unique element
of $\cal C{\nu}_{a}$ in $\omega_a$.

\subsubsection{Itinerary information}\label{sss.combinatorial}

Each element $[z]\in\cal Ck$ has successive sets of good
parameters
$\Gamma_{[z]}^{(k)}\subset\dots\subset\Gamma_{[z]}^{(\nu)}=\Omega_{[z]}$
and corresponding partitions $\cal P{k}_{[z]}, \dots, \cal
P{\nu}_{[z]}$. They are defined in essentially the same way as in
dimension one, as we shall explain in a while. We let $\cal
Pk_{*}$ denote the corresponding induced partition of
$\Omega_{*}^{(k)}$: to each $ \gamma\in\cal Pk_{*} $ is implicitly
associated a critical point $ [z]\in\cal Ck $ with $\gamma\in\cal
Pk_{[z]}$. We always assume that $\gamma$ intersects
$\Gamma^{(k)}$ in at least one point. Otherwise we just delete
$\gamma$: obviously, this has no effect whatsoever on the measure
estimates. Each $\gamma\in\cal Pk_{*}$ has associated
combinatorial information which we call the \emph{itinerary} of
$\gamma$. This consists of:
\begin{itemize}
\item A sequence of \emph{escape times}
\[
 \nu=\eta_{0}<\eta_{1}<\dots<\eta_{s} < k \qquad s=s(\gamma)\geq 0.
\]
\item  Between any two escape times $\eta_{i-1}$ and $\eta_{i}$
 (and between $ \eta_{s} $ and $ k $)
  there is a  sequence  of \emph{essential returns}
\[
 \eta_{i-1}<\nu_{1}<\dots<\nu_{t}<\eta_{i} \qquad t=t(\gamma, i)\geq 0.
\]
\item Between any two essential returns $\nu_{j-1}$ and $\nu_{j}$
 (and between $\nu_{t}$ and $\eta_{i}$) there is a sequence of
 \emph{inessential returns}
\[
 \nu_{j-1}<\mu_{1}<\dots<\mu_{u}<\nu_{j}
 \qquad u=u(\gamma, i, j)\geq 0.
\]
\end{itemize}
Any of these sequences may be empty, except for the first one
because the construction ensures that $\nu$ is always an escape
time.
Any iterate $j$ after an escape time and before the subsequent
return, including the escape time itself, is called an {\em escape
situation\/}. The corresponding image curve
$\gamma_j=\{z_j^{(k)}(a):a\in\gamma\}$ is admissible and long.
Escape times and the essential and inessential return times are
\emph{free returns}. Any returns to $ \Delta$ occurring during
binding periods associated to a previous return are called
\emph{bound returns}. Binding periods for all the returns may be
chosen constant on the interval $\gamma$.

Associated to each free return $\nu$ is a positive integer $|r|$
that we call the \emph{return depth}. This corresponds to the
position of $\gamma_\nu$ relative to the partition $\mathcal
I^{\zeta}=\{I^{\zeta}_{r,m}\}$ as we shall see in the completion
of the inductive
step in Section \ref{ss.proceeding}.
By convention the return depth is zero at escape times. We
let
\[
\cal R{k-1}_{[z]} : \Omega_{[z]}\to \mathbb N
 \qquad\text{ and }\qquad
\cal E{k-1}_{[z]} : \Omega_{[z]} \to \mathbb N.
\]
be the functions which associate to each $a\in\gamma$,
respectively, the sum of all free (essential and inessential)
return depths and the sum of the essential return depths, both
for returns $\nu\le k-1$. These functions are constant on
partition elements, so they naturally induce functions
$$
\cal R{k-1}: \cal Pk_{*}\to \mathbb N
 \qquad\text{ and }\qquad
\cal E{k-1} : \cal Pk_{*}\to \mathbb N.
$$

\subsubsection{Phase and parameter derivatives}\label{sss.phase}

For each $\gamma\in\cal Pk_{*}$ and associated critical point the
velocity $D_{a} z^{(k)}_k(a)$ of the curve $\gamma\ni a\mapsto
z^{(k)}_{k}(a)$ is uniformly comparable, in argument and magnitude,
to the image of the
most expanded vector $ f^{(k)}(z^{(k)}_0(a)) $ under the
differential $D\Phi_a^k(z_0^{(k)}(a))$.
Using also bounded distortion in phase space \BD{k}, we get
a uniform constant $D>0$ such that for every free iterate $k$ the
curve
\[
\gamma_{k}=\{z_{k}^{(k)}(a): a\in\gamma\}= \{\Phi_{a}^{k+1}
(z^{(k)}(a)) : a\in\gamma\}
\]
is admissible and satisfies
\begin{equation}\label{eq.parameterdistortion}
\frac 1 D \frac{|\tilde\gamma_{k}|}{|\gamma_{k}|}
      \leq \frac{|\tilde\gamma|}{|\gamma|}
      \leq D \frac{|\tilde\gamma_{k}|}{|\gamma_{k}|}
\end{equation}
for any subinterval $\tilde\gamma\subset\gamma$.
See \cite{BenCar91}*{Lemmas~8.1, 8.4} and
\cite{MorVia93}*{Lemmas~11.3, 11.5,
11.6} for proofs of these properties. An important
ingredient (cf. \cite{MorVia93}*{Lemma~11.2} or \cite{Via93}*{Lemma~9.2})
is to prove that critical points vary slowly with the parameter
$a$:
$$
\|D_a z^{(k)}(a)\| \le b^{1/20} \ll 1.
$$
In particular, since \( \gamma \) is connected, this gives
\begin{equation}\label{slow}
|z^{(k)}(a)-z^{(k)}(\tilde a)| \leq b^{1/20} |a-\tilde a|
\end{equation}
for all \( a, \tilde a \) belonging to the same element \( \gamma
\in \cal Pk_{[z]} \).
\subsubsection{Existence of binding points}\label{sss.existenceofbindingpoint}

For each $\gamma\in\cal Pk_{*}$ such that $k$ is a free return
there is $\tilde a\in\gamma\cap\Gamma^{(k)}$ and $\zeta^{(k)}\in \cal
Ck_{\tilde a} $ a suitable binding point for all $z^{(k)}_{k}(a),
a\in\gamma $ satisfying $(*)_k$.
 By suitable we mean that
$z_{k}^{(k)}(a)$ and $\zeta^{(k)}(\tilde a)$
are in tangential position for all
$a\in\gamma$ satisfying the recurrence condition at time $k$. This
corresponds to the condition in Section~\ref{sss.binding}.

Note that a critical point for some \emph{fixed} parameter $\tilde
a$ in the intersection \( \gamma\cap\Gamma^{(k)} \) is used as the
binding critical point for all  $a\in\gamma$, we do not need the
continuation \( \zeta(a) \) of $\zeta(\tilde a)$ to be good for
all the parameters in $\gamma$. This is useful when dealing with
the exclusions in Section~\ref{ss.interaction}: we only need to
remove $\gamma$ if {\em all\/} its parameters have anyhow already
been excluded from the set of good parameters. The reason this is
possible is that the interval $\gamma$ is quite small,
$|\gamma|\le e^{-\kappa_0 k}$ for some constant $\kappa_0$ related
to $\kappa$, and critical points vary slowly with the parameter,
see \eqref{slow}, while $(*)_k$ implies $|\zeta^{(k)}(\tilde a) -
z_{k}^{(k)}(a)| \geq \ega{k}\gg e^{-\kappa_0 k}$. See
\cite{MorVia93}*{pp 65-66} or the second Remark in
\cite{Via93}*{\S~9} for explicit estimates.

Moreover, if $p$ represents the binding period associated to the
return $\nu$, then we have
\begin{equation}\label{eq.bindinggrowth}
\frac{|\gamma_{\nu+p+1}|}{|\gamma_{\nu}|}
 \geq e^{\kappa_3 r} \ge e^{\kappa_4 p} \gg 1 ,
\end{equation}
where $r$ is the return depth, and the constants $\kappa_3$,
$\kappa_4$ depend only on $\kappa$. Indeed, this follows from the
corresponding statement in phase space
(\ref{eq.growthinbindingperiods}) and the property \ref{sss.phase}
that phase and parameter derivatives are uniformly comparable.

\begcom

\subsubsection{On the choice of the critical sets}\label{sss.afford2}

The properties in Sections~\ref{sss.globally} and
\ref{sss.location} mean that a tangency point belongs to the
corresponding critical set only if some ancestor of it has a
suitable continuation in parameter space. More precisely, a
critical point $\zeta^{(k)}(\tilde a)$ is included in the critical
set $\cal Ck_{\tilde a}$ at time $k$ only if there exists
$z^{(k-1)}(\tilde a)\in\cal C{k-1}_{\tilde a}$ such that
\begin{enumerate}
\item $\zeta^{(k)}$ has a continuation as a tangency of order $k$
over the entire element $\gamma$ of $\cal Pk_{[z]}$ that contains
$\tilde a$;

\item $z^{(k-1)}(a)$ is an ancestor to $\zeta^{(k)}(a)$ for all
$a\in\gamma$;

\item $k$ is an escape situation for $z^{(k-1)}$, so that
the image of $\gamma\ni a\mapsto z_k^{(k-1)}(a)$ is an admissible
curve with length $\ge \delta/10$.
\end{enumerate}

The reasons we can afford this additional requirement are, on the
one hand, that the recurrence condition $(*)_{k-1}$ implies that
most iterates are escape situations and, on the other hand, that
the binding construction in Section~\ref{sss.afford1} carries on
uniformly over the whole parameter interval. See
\cite{MorVia93}*{\S~11} and the first Remark in
\cite{Via93}*{\S~9} for the technical verifications.

\endcom

\subsection{The parameter space at time \protect\( n \protect\)}
\label{ss.proceeding}

We now explain how the parameter exclusions are determined and how the
parameter space and the combinatorial structure are ``updated''. Part
of this description involves explaining the way that this structure
is updated to take into account the ``new'' critical points, recall
Section \ref{new}, as well as the refinements of ``old'' critical
points, recall Section \ref{old}. We start with the latter.

By definition every interval $\Omega_{[z]}$, with $[z]\in\cal
C{n}$, also belongs to the extended parameter space at the next
iterate $n+1$.
We define a refinement $\hcal P{n}_{[z]}$ of \( \mathcal
P^{(n)}_{[z]} \) based on the dynamics up to time \( n \) and we
update the itinerary information to time \( n \). Based on this
information we decide to exclude some elements of  $\hcal
P{n}_{[z]}$ and then define
\begin{itemize}
\item the set $\Gamma_{[z]}^{(n+1)}$ to be the union of the remaining
elements, and
\item the partition $\cal P{n+1}_{[z]}$ to be simply $\hcal P{n}_{[z]}$
restricted to $\Gamma^{(n+1)}_{[z]}$.
\end{itemize}
The corresponding critical point function $z^{(n)}(a)$ is replaced
by an improvement $z^{(n+1)}(a)$ of order $n+1$, defined on \(
\Gamma^{(n+1)}_{[z]} \). We proceed to explain these steps in
detail.

\subsubsection{Critical neighborhood and
partitions}\label{sss.partitions}

We begin with defining some partition in the dynamical space.
It is no restriction to let $r_\delta=|\log\delta|$ be an integer.
For every integer $r\geq r_\delta$ let
\[
I_{\pm r}=\{z = (x, y)\in Q : \pm x \in (e^{-r-1}, e^{-r}]\}.
\]
 Now let each $ I_{r} $ be further subdivided into $ r^{2}
$ vertical strips $I_{r, m}$ of equal width. This defines a
partition
\[
\mathcal I= \{ I_{r, m}: |r| \geq r_\delta \text{ and } m \in [1,
r^{2}] \}
\]
of $\Delta$ (disregarding $ \{x=0\} $).
Given $I_{r,m}$ we denote $\hat I_{r, m}= I^{\ell}_{r, m}\cup
I_{r,m}\cup I^{\rho}_{r,m}$ where $I^{\ell}_{r, m}$ and
$I^{\rho}_{r, m}$ are the left and right elements adjacent to
$I_{r,m}$. If $I_{r,m}$ happens to be an extreme element of
$\mathcal I$ we just let $I^{\ell}_{r, m}$ or $I^{\rho}_{r, m}$
denote an adjacent interval of length $\delta/10$.
Given a point $\zeta$ we define an analogous partition $\mathcal
I^{\zeta}$ centered at $\zeta$ simply by translating horizontally
$\mathcal I$. Let $I^\zeta_{r}$, $I^\zeta_{r,m}$, $\hat I_{r, m}^\zeta$, be
the corresponding partition elements.

\subsubsection{Partitioning}\label{sss.partitioning}

Let $[z]\in\cal C{n}$ be fixed. For each $\gamma\in \cal
P{n}_{[z]}$ we distinguish two cases.
We call $n$ a non-chopping time in either of the following
situations:
\begin{itemize}
\item[(a)] $\gamma_{n}\cap \Delta$ is empty or contained in an
 outermost partition element of $ \mathcal I^\zeta $;
\item[(b)] $n$ belongs to the binding period associated to some return
time $ \nu<n $ of $ \gamma $.
\end{itemize}
In both situations we simply let $ \gamma\in \hcal P{n}_{[z]}$.

We say that $n$  is a chopping time in the remaining cases, that
is, if $n$ is a free iterate and  $ \gamma_{n} $ intersects $
\Delta $ significantly. In this case we may need to divide $
\gamma_{n} $ into subintervals to guarantee that the distortion
bounds and other properties continue to hold. This depends on the
position of $\gamma_{n}$
relative to the partition ${\mathcal I}^{\zeta}$ of $\Delta$
centered at the binding point $\zeta$ (recall that $\gamma_{n}$
is admissible, and thus ``transverse'' to the partitions).
Indeed, there are two different situations.

If $ \gamma_{n} $ intersects at most two partition elements then
we let $\gamma\in \hcal P{n}_{[z]}$, that is, we do not
subdivide it. In this case we say that $n$ is an inessential
return and add this to the itinerary information already
associated to $\gamma$. If $\gamma_{n}$ does intersect at least
three partition elements then we partition it
\begin{equation}\label{pullback}
\gamma = \gamma^{\ell}\cup \bigcup_{(r, m)}\gamma^{(r, m)}\cup
\gamma^{\rho},
\end{equation}
where (using notations from Section \ref{sss.partitions})
\begin{itemize}
\item each $\gamma_{n}^{(r, m)}=\{z_{n}^{(n)}(a): a\in\gamma^{(r, m)}\}$
crosses $I_{r,m}^\zeta$ and is contained in $\hat I_{r,m}^\zeta$;
\item $ \gamma_{n}^{\ell} $ and $ \gamma_{n}^{\rho} $ are either empty or
 components of $\gamma_{n}\setminus \Delta$
 with length $\ge \delta/10$.
\end{itemize}
If the connected components of $\gamma_{n}\setminus \Delta$ have
length $<\delta/10$ we just glue them to the adjacent interval of
the form $\gamma_{n}^{(r, m)}$.

By definition the resulting subintervals of $\gamma$ are elements of
$\hcal P{n}_{z}$.
 The intervals $ \gamma^{\ell}$ , $ \gamma^{\rho} $ are called \emph{escape components}
and are said to have an escape at time $ n $. All the other
intervals are said to have an essential return at time $ n $ and
the corresponding values of $ |r| $ are the associated essential
return depths. By convention, escaping components have return
depth zero.

\subsubsection{Parameter exclusions}

We now have itinerary information up to time $n$ and in
particular return depths are defined and so are the functions
$\cal R{n}_{[z]}$ and $\cal E{n}_{[z]}$ on
 $\hcal P{n}_{[z]}$.
We fix a small constant $\tau >0 $ and then let
\begin{equation}\label{eq.exclusions}
E^{(n)}_{[z]}
 = \bigcup \{\gamma\in\hcal P{n}_{[z]}:
 \cal E{n}_{[z]}(\gamma) > \tau n \}.
\end{equation}
We claim that if $\tau$ is small with respect to $\alpha$ then all
$\gamma\in\hcal P{n}_{[z]}$ that are {\em not\/} in
$E^{(n)}_{[z]}$ satisfy the recurrence condition $(*)_{n}$.
The proof is postponed to Section~\ref{sss.verifying}. We define
\begin{equation}\label{eq.newgoodset1}
  \Gamma_{[z]}^{(n+1)} = \Gamma_{[z]}^{(n)} \setminus
  \left(E_{[z]}^{(n)} \cup \{\gamma\in \hcal P{n}_{[z]}:
                        \gamma\cap
                        \Gamma^{(n)}=\emptyset\}\right).
\end{equation}
Moreover, we define $\cal P{n+1}_{[z]}$ as the restriction of $\hcal
P{n}$ to $\Gamma^{(n+1)}_{[z]}$.

Observe that we also remove from
the parameter space of $[z]$ intervals that have already been
completely deleted because of other critical points, even if they
may look like good parameters for the point $[z]$. This deals with
the interaction between distinct critical orbits discussed in
Section~\ref{ss.interaction}. Let us stress once more that these
additional exclusions have no effect at all on the new set
of parameters good for all critical points:
\begin{equation}\label{eq.newgoodset2}
  \Gamma^{(n+1)}= \Gamma^{(n)}\setminus \!\!\!
  \bigcup_{[z]\in\cal C{n}} \!\!\! E_{[z]}^{(n)}.
\end{equation}
For parameters $a\in\Gamma_{[z]}^{(n+1)}$ we replace $z^{(n)}(a)$
by its refinement $z^{(n+1)}(a)$ as discussed in Sections
\ref{ss.new critical points}  and
\ref{sss.newcritical}. This completes the construction in this
case.

\subsubsection{Verifying the recurrence condition}\label{sss.verifying}

We explain why all $[z]=z^{(n)}\in\cal C{n}$ satisfy the condition
$ (*)_{n}$ for all $a\in\Gamma_{[z]}^{(n+1)}$. We begin with the
remark that if $\nu$ is a return for $[z]$ then, by construction,
the return depth $|r_\nu| \approx |\log d(z_{\nu})|$. In
particular the sum in $(*)_{n}$ is
 $\le 2\cal R{n}([z])$. We claim moreover that the sum of all free
return depths is bounded by a multiple of the sum of essential
return depths: there exists a uniform constant $C>0$ such that
\begin{equation}\label{eq.depths}
 \cal R{n}([z]) \leq C \cal E{n}([z]).
\end{equation}
Assuming this statement, and keeping (\ref{eq.exclusions}) and
(\ref{eq.newgoodset2}) in mind, we get that for all the parameters
$a\in\Gamma_{[z]}^{(n+1)}$
$$
\sum_{\substack{\text{free}\\\text{returns} \\ \nu \leq n}}
\!\!\! \log \frac{1} {d(z_{\nu}^{(n)})}
 \leq 2 \cal R{n}([z])
 \leq 2 C \cal E{n-1}([z])
 \leq 2 C \tau.
$$
The conclusion follows choosing $\tau < \alpha/(2C)$.

We are left to prove (\ref{eq.depths}). Let $\mu_1<\cdots<\mu_s$
be the inessential returns in between consecutive essential
returns $\nu_i<\nu_{i+1}$. Also let $r_i$ be the return depth
associated to $\nu_i$ and $\rho_j$ be the return depth associated
to each $\mu_j$, $1\le j\le s$. Property (\ref{eq.bindinggrowth})
says that the iterates $|\gamma_j|$ are expanded over the complete
binding period associated to any free return. Due to the
hyperbolic behavior of our maps outside the critical region, we
know that these curves are not contracted during free periods.
This gives that
$$
|\gamma_{\nu_{i+1}}|
 \ge e^{\kappa_3(r_i+\rho_1+\cdots+\rho_s)}|\gamma_{\nu_{i}}|.
$$
Clearly, $|\gamma_{\nu_{i+1}}|\le 2$. On the other hand,
$$
|\gamma_{\nu_{i}}| \ge \const e^{-r_i}r_i^{-2} \ge 2 e^{-2 r_i}
$$
because $\nu_i$ is an essential return. Putting these two
estimates together we find
$$
e^{\kappa_3(r_i+\rho_1+\cdots+\rho_s)-2 r_i} \le 1
\quad\Rightarrow\quad
 \rho_1+\cdots+\rho_s \le (2/\kappa_3) r_i\,.
$$
Adding these inequalities for every essential return $\nu_i$, we
get $$\cal R{n}([z]) \leq C \cal E{n}([z])$$ with $C=2/\kappa_3$.

\subsubsection{New critical points}\label{sss.newcritical}

Finally, we must include in the construction new critical points
of order \( n+1 \).   First of all we
``upgrade''  the old critical points \( [z]=z^{(n)} \) as described in
Section \ref{old}. These are easily seen to satisfy the inductive
assumptions stated in Sections \ref{sss.globally} and
\ref{sss.location}. Then we add ``really new'' critical
points $[\zeta]=\zeta^{(n+1)}$ as mentioned in Section
\ref{new}. To ensure that the inductive assumptions continue to hold
for these points we proceed as follows.
For every critical
point $[z]=z^{(n)}$ of order $n$ and every $\gamma$ in the
corresponding partition $\cal P{n+1}_{[z]}$ such that $n$ is an escape
situation, as defined in the previous section, we introduce the
points $\zeta^{(n+1)}$ such that
\begin{itemize}
\item $\zeta^{(n+1)}$ is defined as a critical continuation over the whole
$\gamma$;
\item $\zeta^{(n+1)}(a)$ is contained in
      $\Phi_a^{\theta (n+1)}(W(a))\setminus\Phi_a^{\theta n}(W(a))$
      for all $a\in\gamma$;
\item $z^{(n)}(a)$ is ancestor to $\zeta^{(n+1)}(a)$ for all
$a\in\gamma$.
\end{itemize}
Essentially, these $\zeta^{(n+1)}$ are the additional elements of
the critical set at time $n+1$. However, there is the possibility
that two or more critical points \( z^{i,(n)} \), defined on
intervals $\gamma_i$\,, generate by this procedure critical
functions $\zeta^{i,(n+1)}$ which turn out to coincide at some
parameters:
\begin{equation}\label{redundance}
\zeta^{i,(n+1)}(a) = \zeta^{j,(n+1)}(a) \quad\text{for some (and
hence all) $a\in\gamma_i\cap\gamma_j$\,.}
\end{equation}
If we were to consider all these $\zeta^{i,(n+1)}(a)$ as different
critical points, the counting argument to prove \eqref{meas1} that
we give in the next paragraph would not be valid. Instead, we
begin by (almost) eliminating redundancy as follows. From any
family of critical functions $\zeta^{i,(n+1)}$ as in
\eqref{redundance} we extract a minimal subfamily
$\zeta^{i_j,(n+1)}$ such that the union of their domains
$\gamma_{i_j}$ coincides with the union of all $\gamma_i$\,. We
retain these $\zeta^{i_j,(n+1)}$ but eliminate all the other
$\zeta^{k,(n+1)}$ as they are clearly redundant. The key, if quite
easy observation is that by minimality a parameter $a$ belongs to
not more than two of these domains $\gamma_{i_j}$.

The critical functions $\zeta^{(n+1)}$ obtained in this way, after
the redundancy elimination we just described, are the remaining
elements $[\zeta]$ of $\cal C{n+1}$. For each one of them we set
$$
\Omega_{[\zeta]}=\Gamma_{[\zeta]}^{(n+1)}=\gamma
\quad\text{and}\quad\cal P{n+1}_{[\zeta]}=\{\Omega_{[\zeta]}\}.
$$
We think of these critical points as being ``born'' at time \( n
\). Thus the iterate $\nu=n$ is the first escape time for each
$[\zeta]$; apart from this the combinatorics of $[\zeta]$ is blank
and there are no exclusions corresponding to these points at this
time. This procedure \emph{defines} the new critical set $\cal
C{n+1}$ and, in particular, makes precise the meaning of the
symbol $[z]$ at the next stage of the construction. By the
condition in Section~\ref{sss.location}, there are at most
$$
\frac{2|\Phi_a^{\theta (n+1)}(W(a))
\setminus\Phi_a^{\theta(n)}(W(a))|} {2\rho^{\theta (n+1)}}
$$
of these new critical points whose domains $\Omega_{[\zeta]}$
contain a given $a\in\Omega$ (here we have $\nu=n$). The factor
$2$ in the numerator accounts for the fact that a given point may
represent two ``different'' critical points, but not more, as
discussed in the previous paragraph. In other words,
\begin{equation*}
    \begin{aligned}
{}&
\#\{[w]\in\cal C{n+1}: a\in\Omega_{[w]}\} \le
 \\
 & \qquad\qquad \le \#\{[z]\in\cal C{n}:
 a\in\Omega_{[z]}\} +
    \frac{|\Phi_a^{\theta (n+1)}(W(a))
    \setminus\Phi_a^{\theta(n)}(W(a))|}
    {\rho^{\theta (n+1)}}\,.
\end{aligned}
\end{equation*}
Now a simple induction argument yields the bound in \eqref{meas1}
$$
 \#\{[w]\in\cal C{n+1}: a\in\Omega_{[w]}\}
 \le\frac{|\Phi_a^{\theta (n+1)}(W(a))|}{\rho^{\theta (n+1)}}
 \le (5/\rho)^{\theta (n+1)}.
$$
The $n$'th step of the construction is complete.

\section{The probabilistic argument}\label{s.exclusionestimates}

It remains to show that the set $\Gamma^{*}=\cap_{n}\Gamma^{(n)}$
has positive Lebesgue measure. For each critical point we use a
large deviations argument similar to the one-dimensional proof to get
the estimate as in \eqref{meas2}.
As discussed in Section \ref{sss.parexc}, we then sum the exclusions
associated to each critical point, using the
bound on the number of critical points in \eqref{meas1}.
Thus, most of this section deals with the exclusions associated to a
single critical point  $[z]\in\cal Cn$ for some $n \ge N$. The issue
of the multiplicity of critical point is taken care of by the formalism.
For simplicity we omit the subscript \( [z] \) where this does not
give rise to confusion.

We split the argument into 4 sections. The first step is a useful
re-organization of the combinatorial structure on each
\(\Omega_{[z]}\). The reason this is necessary is that our
combinatorial data keeps track of the critical orbits itineraries
in between escape returns (escape situations which coincide with
chopping times) but not beyond. At escape returns the dynamics
starts ``afresh", in the sense that subsequent itineraries are
very much independent of the previous behavior. This is a key
feature of escape returns (the system ``escapes its past") and,
together with the fact that such returns are fairly frequent
(large waiting time exponentially improbable), a crucial
ingredient in the whole exclusion argument. On the other hand it
means that, due to the possibility of many intermediate escapes, the
same combinatorial data may correspond to several different
partition intervals, even with unbounded multiplicity. The purpose
of the re-organization we carry out in Section
\ref{ss.combinatorics} is to decompose the whole collection of
elements of all the partitions into a number of ``blocks'' on each
of which we do have bounded multiplicity of the combinatorics. The
strategy is to restrict our attention to the subintervals of an
escaping component only up to the following escape time associated
to each subinterval.

Focussing on each one of these blocks, we show in Section
\ref{ss.metric} that intervals are exponentially small in terms of
the total sum of their return depths. Then in Section
\ref{ss.combinatorial} we develop a counting argument to estimate
the cardinality of the set of intervals whose return depths sum up
to some given value. We show that this bound is exponentially
increasing in the sum of the return depths, but with an
exponential rate slower than that used to estimate the size of the
intervals. Combining these two estimates immediately implies a
bound on the \emph{average} recurrence for points in a single
block. In Corollary~\ref{c.mainestimate} we then show how to
``sum'' the contributions of each block to get an estimate of the
overall average recurrence over all points of \( \Gamma^{(n)} \).
Finally, a large deviation argument implies the required estimate
for the proportion of excluded parameters.

\subsection{Combinatorics renormalization}
\label{ss.combinatorics} By construction, for each $\gamma\in\hcal
Pn_{[z]}$ we have a sequence
$\nu=\eta_{0}<\eta_{1}<\dots<\eta_{s}\leq n$, $s\geq 0$ of escape
times and for each $0\leq i\leq s$ there exists an ancestor
$z^{(\eta_{i})}\in\cal C{\eta_{i}}$ and an interval
$\gamma^{(\eta_{i})}$ with
$\gamma\subset\gamma^{(\eta_{i})}\subset\Omega_{[z]}$ and which is
an escape component for $z^{(\eta_{i})}$. In particular,
$z^{(\eta_{i})}$ admits a continuation to the whole
$\gamma^{(\eta_{i})}$. Because $s$ may depend on $\gamma$, it is
convenient to extend the definition of $\gamma^{(\eta_{i})}$ to
all $ 0\leq i\leq n$ and we do this by letting
$\gamma^{(\eta_{i})}=\gamma$ for $s+1\leq i\leq n $. Then we
consider the disjoint union
\[
\cal Qi = \cal Qi_{[z]}
 = \coprod\{\gamma^{(\eta_{i})} : {\gamma\in\hcal Pn_{[z]}}\}
\]
Clearly, $\cal Q0 = \{\Omega_{[z]}\}$ and $\cal Qn
=\hcal Pn_{[z]}$, which is a partition of $\Gamma_{[z]}^{(n)}$.
Given $\gamma^{(\eta_{i})}\in\cal Qi$ and
$\gamma^{(\eta_{j})}\in\cal Qj$ with $i<j$, we say that
$\gamma^{(\eta_{j})}$ is a \emph{descendant} of
$\gamma^{(\eta_{i})}$ if
$\gamma^{(\eta_{j})}\subset\gamma^{(\eta_{i})}$ and
$z^{(\eta_{i})}$ is ancestor to $ z^{(\eta_{j})}$. For $0\leq i
\leq n-1$ and $\gamma\in \cal Qi$ we let
\[
\cal Q{i+1}(\gamma) = \{\gamma'\in\cal Q{i+1} : \gamma' \text{ is
a descendant of } \gamma \}.
\]
The itineraries of all intervals in $ \cal Q{i+1}(\gamma) $
clearly coincide up to time $ \eta_{i} $. Then we may define
functions
 $\Delta\cal Ei_{\gamma}: \cal Q{i+1}(\gamma) \to \mathbb N$ where
\[
\Delta\cal Ei_{\gamma}(\gamma')
 = \cal E{\eta_{i+1}}(\gamma')-\cal E{\eta_{i}}(\gamma)
\]
is the sum of all essential return depths associated to the
itinerary $\gamma'\in\cal Q{i+1}(\gamma)$ between the escape times
$\eta_{i}$ and $\eta_{i+1}$. Finally we let
\[
 \cal Q{i+1}(\gamma, R)
  = \{\gamma'\in \cal Q{i+1}(\gamma): \Delta\cal Ei_{\gamma} (\gamma')= R\}.
\]

\subsection{Metric bounds}\label{ss.metric}

Let $\bar\kappa=\kappa_3/5$, where $\kappa_3$ is the constant in
\eqref{eq.bindinggrowth}. Recall that $\kappa_3$ is independent of
$\delta$ and $b$.

\begin{lemma}\label{l.metric}
For all $[z]\in\cal C{n}$,  $0\leq i \leq n-1$, $\gamma\in \cal
Qi$, $R\geq 0$, and $\gamma'\in \cal Qi(\gamma,R)$ we have
$$
|\gamma'|\leq e^{-3\bar\kappa R} |\gamma|.
$$
\end{lemma}

\begin{proof}
By construction there are nested intervals
$\gamma'\subset\gamma^{(\nu_{t})}\subset
 \dots\subset\gamma^{(\nu_{1})}\subset\gamma^{(\nu_{0})}
 =\gamma$
such that $\nu_0$ is an escape time for
$\gamma=\gamma^{(\nu_{0})}$ and for each $j=1, \dots, t$ the
interval $\gamma^{(\nu_{j-1})}$ has an essential return at time
$\nu_{j}$ which is when the interval $\gamma^{(\nu_{j})}$ is
created as a consequence of chopping. Write
\begin{equation}\label{eq.telescope0}
\frac{|\gamma'|}{|\gamma|}=\frac{|\gamma^{(\nu_{1})}|}
                                     {|\gamma^{(\nu_{0})}|}
                                \frac{|\gamma^{(\nu_{2})}|}
                                     {|\gamma^{(\nu_{1})}|}
                          \dots \frac{|\gamma^{(\nu_{t})}|}
                                     {|\gamma^{(\nu_{t-1})}|}
                                \frac{|\gamma'|}
                                     {|\gamma^{(\nu_{t})}|}.
\end{equation}
The last factor has the trivial bound
$|\gamma'|/|\gamma^{(\nu_{t})}|\leq 1$. For the middle factors we
use

\begin{lemma}\label{lemma4}
For all $j= 1, \dots, t-1$ we have
\[
\frac{|\gamma^{(\nu_{j+1})}|}{|\gamma^{(\nu_{j})}|}
 \leq e^{-r_{j+1}+ (1-3\bar\kappa)r_{j}}\,.
\]
\end{lemma}

\begin{proof}
By the bounded distortion property \eqref{eq.parameterdistortion},
\begin{equation}\label{eq.telescope1}
\frac{|\gamma^{(\nu_{j+1})}|}{|\gamma^{(\nu_{j})}|}
 \leq D \frac{|\gamma^{(\nu_{j+1})}_{\nu_{j}+p_{j}+1}|}
             {|\gamma^{(\nu_{j})}_{\nu_{j}+p_{j}+1}|}\,.
\end{equation}
For each of these essential returns \eqref{eq.bindinggrowth}
gives,
\[
|\gamma^{(\nu_{j})}_{\nu_{j}+p_{j}+1}|
 \geq e^{5\bar\kappa r_j}|\gamma^{(\nu_{j})}_{\nu_{j}}|
 \geq \const r_j^{-2} e^{(5\bar\kappa - 1) r_{j}}
 \geq e^{(4\bar\kappa - 1) r_{j}}.
\]
We used here that $r^{2}$ is much smaller than $e^{\bar\kappa r}$
for $r\ge r_\delta\gg 1$. To get an upper bound for the numerator
in \eqref{eq.telescope1} we use that $\gamma^{(\nu_{j+1})}_{k}$
remains outside the critical region and is an admissible curve
between time $\nu_{j}+p_{j}$ and time $\nu_{j+1}$. So, during this
period its length can not decrease:
\[
|\gamma^{(\nu_{j+1})}_{\nu_{j}+p_{j}+1}|
 \leq |\gamma^{(\nu_{j+1})}_{\nu_{j+1}}|
 \leq e^{-r_{j+1}}.
\]
Replacing these two bounds in \eqref{eq.telescope1}, we find
(using that $r_\delta$ is large)
$$
\frac{|\gamma^{(\nu_{j+1})}|}{|\gamma^{(\nu_{j})}|}
 \le D e^{-r_{j+1}+(1-4\bar\kappa)r_j}
 \le e^{-r_{j+1}+(1-3\bar\kappa)r_j}
$$
as claimed in Lemma~\ref{lemma4}.
\end{proof}

A similar argument applies to the first factor ($j=0$) of
\eqref{eq.telescope0}. The length of
$\gamma^{(\nu_{1})}_{\nu_{1}}$ is bounded by $e^{-r_1}$, by
construction. Moreover, the escaping component
$\gamma^{(\nu_{0})}_{\nu_{0}}$ has length $\ge\delta/10$. Since
this component is adjacent to $\Delta$ and all the iterates from
time $\nu_0$ to time $\nu_1$ take place in the hyperbolic region
$\Delta^c$, a simple hyperbolicity argument gives
$|\gamma^{(\nu_{0})}_{\nu_{1}}|\ge\delta^{9/10}$. It is no
restriction to suppose $5\bar\kappa=\kappa_3<1/10$. Thus, we get
\begin{equation}\label{eq.telescope2}
\frac{|\gamma^{(\nu_{1})}|}
     {|\gamma^{(\nu_{0})}|}
 \leq D \, \frac{|\gamma^{(\nu_{1})}_{\nu_{1}}|}
             {|\gamma^{(\nu_{0})}_{\nu_{1}}|}
 \leq D \, e^{-r_1} \delta^{-1+5\bar\kappa}
 \leq e^{-r_1} \delta^{-1+3\bar\kappa}.
\end{equation}
Replacing these bounds in \eqref{eq.telescope0} we find
\begin{align*}
\frac{|\gamma'|}{|\gamma|}
  & \leq \exp\big(- r_1 - (1-3\bar\kappa)\log\delta + \sum_{j=1}^{t-1}
                  - r_{j+1} + \sum_{j=1}^{t-1} (1-3\bar\kappa)r_j \big)
  \\
  & = \exp\big(- (1-3\bar\kappa)\log\delta - r_t - 3\bar\kappa\sum_{j=1}^{t-1} r_j \big)
  \leq \exp\big(-3\bar\kappa \sum_{j=1}^{t} r_{j}\big).
\end{align*}
In the second inequality we have used $r_t\ge|\log\delta|$. The
term on the right is $e^{-3\bar\kappa R}$, so this completes the
proof of Lemma~\ref{l.metric}.
\end{proof}

\subsection{Combinatorial bounds}\label{ss.combinatorial}

\begin{lemma}\label{l.combinatorial}
For all $[z]\in\cal C{n}$,  $0\leq i \leq n-1$, $\gamma\in \cal
Qi$ and $ R\geq 0 $ we have
\[
\# \cal Q{i+1}(\gamma, R) \leq  e^{\bar\kappa R}
\]
\end{lemma}

\begin{proof}
Given $\gamma'\in\cal Q{i+1}(\gamma, R)$ let
$\eta_i=\nu_0<\nu_1<\cdots<\eta_t<\eta_{i+1}$ be the corresponding
sequence of essential returns between the consecutive escape
situations. To each $\nu_i$ the chopping procedure in
Section~\ref{sss.partitioning} assigns a pair of integers
$(r_i,m_i)$ with $|r_i|\geq r_\delta$ and $1\le m_i\le r_i^2$. By
definition $|r_{1}|+\dots+|r_{t}|=R$. The sequence $(r_i,m_i)$
determines $\gamma'$ completely, except that the next escape
situation $\eta_{i+1}$ may be generating two escape components,
which thus share the same sequence. So, apart from the harmless
factor $2$, the cardinality of $\cal Q{i+1}(\gamma, R)$ is bounded
by the number of sequences $(r_i,m_i)$, $t\ge 1$, with
$|r_{1}|+\dots+|r_{t}|=R$ and $|r_i|\geq r_\delta$ and $1\le
m_i\le r_i^2$.

We begin by estimating the number of integer solutions to
\begin{equation}\label{eq.stirling}
r_{1} + \dots + r_{t}=R\quad\text{with}\quad r_i \geq r_\delta.
\end{equation}
This corresponds to the number of ways of partitioning $ R $
objects into $t$ disjoint subsets, which is well known to be
bounded above by $(R+t)!/R!t!$. Using Stirling's approximation
formula we have
\begin{equation*}
\frac{(R+t)!}{R! t!}
 \leq \const \frac{(R+t)^{R+t}}{R^{R}t^{t}}
 = \const \left[\left(1+\frac tR\right)^{1+\frac tR} \left(\frac Rt\right)^{\frac tR}\right]^R.
\end{equation*}
Recalling the fact that $R\ge t\, r_\delta$ and
$r_\delta\to\infty$ when $\delta\to 0$, we see that both factors
in the last term tend to 1 when $\delta$ tends to zero. Therefore,
choosing $\delta$ small enough we ensure that the number of
solutions of \eqref{eq.stirling} is less than $e^{\bar\kappa
R/4}$.

Now to complete the proof we need to  sum over all values of $t$
and we also need to take into account the sign of each $r_i$ and
the variation of $m_i$ from $1$ to $r_{i}^{2}$. The latter means
that each fixed sequence $r_i$ corresponds to at most
$\prod_{i=1}^{t}r_{i}^{2}$ partition elements. For fixed $R$ the
product is biggest when the $r_i$ are approximately equal. So this
multiplicity is bounded by $(R/t)^t$. In the range we are
interested in, the function $t\mapsto(R/t)^t$ is monotone
increasing on $t$.  So we may bound it by $r_\delta^{R/r_\delta}$,
which is $\le e^{\bar\kappa R/4}$ if $r_\delta$ is large. In this
way we get
\begin{equation*}
\# \cal Q{i+1}(\gamma,R)
 \leq  2 \sum_{t \leq R/r_{\delta}} 2^{t}  \, e^{\bar\kappa R/2}
 \leq  4 \, 2^{R/r_{\delta}} \, e^{\bar\kappa R/2}
 \leq e^{\bar\kappa R}
\end{equation*}
if $r_\delta$ is large enough.
\end{proof}

\subsection{Average recurrence}\label{ss.average}

From Lemmas~\ref{l.metric} and \ref{l.combinatorial} we
immediately get
\begin{equation}\label{eq.mainestimate}
\sum_{\gamma'\in\cal Q{i+1}(\gamma, R)}|\gamma'|
 \leq e^{-2\bar\kappa R} |\gamma|.
\end{equation}
With the help of this we are going to give an estimate for the
distribution of the recurrence function when the entire itinerary
up to time $n$ is taken into account.

\begin{corollary}\label{c.mainestimate}
For every $[z]\in\cal C{n}$ we have
\[
\int_{\Gamma^{(n)}_{[z]}} e^{\bar\kappa\cal En}
 \leq e^{n/r_\delta}|\Omega_{[z]}|.
\]
\end{corollary}
\begin{proof}
Recall first of all that $\cal Qn=\hcal Pn_{[z]}$ is a partition
of $\Gamma^{(n)}_{[z]} $ and that $\cal En$ is constant on
elements of $\cal Qn$. Thus
$$
\int_{\Gamma^{(n)}_{[z]}} e^{\bar\kappa\cal En}
 = \sum_{\gamma^{(n)}\in\cal Qn} e^{\bar\kappa\cal En} |\gamma^{(n)}|\,.
$$
Moreover $\cal En=\cal E{n-1}+\Delta\cal E{n-1}$ where
 $\cal E{n-1}$ depends only on the element $\gamma^{(n-1)}$ of $\cal Q{n-1}$ containing
 $\gamma^{(n)}$. Therefore
\[
\sum_{\gamma^{(n)}\in\cal Qn} e^{\bar\kappa \cal En}|\gamma^{(n)}|
 \leq \sum_{\gamma^{(n-1)}\in\cal Q{n-1}} e^{\bar\kappa\cal E{n-1}}
 \sum_{\gamma^{(n)}\in\cal Qn} e^{\bar\kappa \Delta\cal En}
 |\gamma^{(n)}|
\]
and iterating the argument
\begin{equation}\label{eq.split}
\begin{split}
\sum_{\gamma^{(n)}\in\cal Qn} e^{\bar\kappa\cal En} |\gamma^{(n)}|
 \leq
\sum_{\gamma^{(1)}\in \cal Q1(\gamma^{(0)})}
e^{\bar\kappa\Delta\cal E0} \dots \sum_{\gamma^{(i+1)}\in \cal
Q{i+1}(\gamma^{(i)})} e^{\bar\kappa\Delta\cal E{i}}\dots
\\
\dots\sum_{\gamma^{(n-1)}\in \cal Q{n-1}(\gamma^{(n-2)})}
e^{\bar\kappa\Delta\cal E{n-2}} \sum_{\gamma^{(n)}\in \cal
Q{n}(\gamma^{(n-1)})} e^{\bar\kappa\Delta\cal E{n-1}}
|\gamma^{(n)}|
\end{split}
 \end{equation}
For each $  0\leq i \leq n-1$ and $ \gamma^{(i)}\in \cal Qi $ we
can write
\[
 \sum_{\gamma^{(i+1)}} e^{\bar\kappa\Delta\cal Ei(\gamma^{(i+1)})} |\gamma^{(i)}|
 = |Q^{i+1}(\gamma^{(i)}, 0)| + \sum_{R\geq r_\delta} e^{\bar\kappa R} |Q^{(i+1)} (\gamma^{(i)}, R)|
\]
where the sum on the left is over all
 $\gamma^{(i+1)}\in\cal Q{i+1}(\gamma^{(i)})$. The relation \eqref{eq.mainestimate} gives
\[
 \sum_{R\geq r_\delta} e^{\bar\kappa R} | Q^{(i+1)} (\gamma^{(i)}, R)|
  \leq \sum_{R\geq r_\delta} e^{-\bar\kappa R} |\gamma^{(i)}|
  \leq 2 e^{-\bar\kappa r_\delta} |\gamma^{(i)}|
\]
and therefore
\begin{equation}\label{eq.plug}
 \sum_{\gamma^{(i+1)}}
 e^{\bar\kappa\Delta\cal Ei(\gamma^{(i+1)})} |\gamma^{(i)}|
 \leq (1+2 e^{-\bar\kappa r_\delta})|\gamma^{(i)}|
 \leq e^{1/r_\delta} |\gamma^{(i)}|
\end{equation}
(the last inequality uses that $\log(1+2 e^{-\bar\kappa r_\delta})
\le 4 e^{-\bar\kappa r_\delta} \le 1/r_\delta$ assuming $r_\delta$
is large enough). Replacing \eqref{eq.plug} in \eqref{eq.split},
successively for all values of $i$, and recalling that
$\gamma^{(0)}=\Gamma^{(0)}_{[z]}=\Omega_{[z]}$, we find
$$
\sum_{\gamma^{(n)}\in\cal Qn} e^{\bar\kappa\cal En} |\gamma^{(n)}|
\leq e^{n/r_\delta}|\Omega_{[z]}|
$$
as claimed.
 \end{proof}

\subsection{Conclusion}\label{ss.conclusion}

Using the Chebyshev inequality and the definition
 \[
E^{(n)}_{[z]}
 = \{\gamma\in\hcal Pn_{[z]}: \cal En \geq \tau(n+1)\}\]
we get (recall that $\hcal Pn_{[z]}$ is a partition of
$\Gamma^{(n)}_{[z]}$)
\[
|E^{(n)}_{[z]}|
 \leq e^{-\tau \bar\kappa(n+1)}\int_{\Gamma^{(n)}_{[z]}} e^{\bar\kappa \cal En}
 \leq e^{\left(\frac{1}{r_\delta}-\bar\kappa\tau\right) n}|\Omega_{[z]}|
 \leq e^{-\bar\kappa\tau n/2}|\Omega_{[z]}|.
\]
Then, using \eqref{meas1},
\[
|\bigcup_{[z]\in\cal Cn} E^{(n)}_{[z]}|
 \leq e^{-\bar\kappa\tau n/2} \sum_{[z]} |\Omega_{[z]}|
 \leq e^{-\bar\kappa\tau n/2}(5/\rho)^{\theta n} |\Omega| \,.
\]
While $\bar\kappa$ and $\tau$ are independent of the perturbation
size $b$, the constant $\theta$ can be made arbitrarily small by
reducing $b$. So, we may suppose that the last expression is $\le
e^{-\bar\kappa\tau n/4} |\Omega|$.

This means that
 $|\Gamma^{(n)}\setminus\Gamma^{(n+1)}| \le e^{-\bar\kappa\tau n/4} |\Omega|$
for all $n$, which implies
\[
 |\Gamma^{(n+1)}|
  \geq \Big(1-\sum_{i=N}^{n}e^{-\bar\kappa\tau i/4}\Big) |\Omega|
\]
and
\[
 |\Gamma^{*}|
  \geq \Big(1-\sum_{n=N}^{\infty}e^{-\bar\kappa\tau n/4}\Big) |\Omega| > 0.
\]

\end{document}